 \documentclass{article}
 \usepackage{amstex}
 \usepackage{amssymb}
 \usepackage{theorem}

\newenvironment{pf}{\paragraph{Proof}}{\par\smallskip}            
\newenvironment{pfof}[1]{\paragraph{Proof of #1}}{\par\smallskip} 
\newtheorem{TEO}{Theorem}[section]
\newtheorem{COR}[TEO]{Corollary}
\newtheorem{PROP}[TEO]{Proposition}

\newtheorem{LEM}[TEO]{Lemma}
\newtheorem{CLA}[TEO]{Claim}
{
\theorembodyfont{\rmfamily}
\newtheorem{DEF}[TEO]{Definition}
\newtheorem{REM}[TEO]{Remark}
\newtheorem{REMS}[TEO]{Remarks}
}

\newcommand\Oh{{\cal O}}
\newcommand\sC{{\cal C}}

\newcommand\sF{{\cal F}}
\newcommand\sG{{\cal G}}
\newcommand\sI{{\cal I}}
\newcommand\sJ{{\cal J}}
\newcommand\sL{{\cal L}}
\newcommand\sN{{\cal N}}

\newcommand\al{\alpha}
\newcommand\be{\beta}
\newcommand\ga{\gamma}
\newcommand\om{\omega}
\newcommand\la{\lambda}
\newcommand\La{\Lambda}
\newcommand\si{\sigma}
\newcommand\Ga{\Gamma}
\newcommand\de{\delta}

\newcommand\De{\Delta}
\newcommand\fie{\varphi}

\newcommand\simby[1]{\buildrel{{\scriptstyle\mathrm{#1}}}\over{\sim}}
\newcommand\numeq{\simby{{num}}}         
\newcommand\lineq{\simby{{lin}}}         

\newcommand\dual{\mathrel{\raise3pt\hbox{$\underline{\mathrm{\thinspace d
\thinspace}}$}}}
\newcommand\1{^{-1}}

\newcommand\Bh{\widehat B}
\newcommand\QED{\ifhmode\unskip\nobreak\fi\quad {\rm Q.E.D.}} 
\newcommand\qed{\ifhmode\unskip\nobreak\fi\quad $\Box$}       
\newcommand\iso{\cong}
\newcommand\into{\hookrightarrow}
\newcommand\onto{\twoheadrightarrow}
\newcommand\bij{\leftrightarrow}
\newcommand\broken{\mathrel{{\relbar\kern-.2pt\rightarrow}}}
\newcommand\aff{\Bbb A}
\newcommand\proj{\Bbb P}
\newcommand\Z{\Bbb Z}
\newcommand\id{\mathrm{id}}
\newcommand\res{\mathrm{res}}

\newcommand\length{\operatorname{length}}
\newcommand\chara{\operatorname{char}}
\newcommand\coker{\operatorname{coker}}

\newcommand\rank{\operatorname{rank}}

\newcommand\Ann{\operatorname{Ann}}
\newcommand\Bs{\operatorname{Bs}}
\renewcommand\div{\operatorname{div}}
\newcommand\Ext{\operatorname{Ext}}
\newcommand\sExt{\operatorname{{\cal E}{\it xt}}}
\newcommand\Hom{\operatorname{Hom}}
\newcommand\End{\operatorname{End}}
\newcommand\sHom{\operatorname{{\cal H}{\it om}}}

\newcommand\Norm{\operatorname{Norm}}

\newcommand\Pico{\operatorname{Pic^0}}
\newcommand\Sing{\operatorname{Sing}}
\newcommand\NonSing{\operatorname{NonSing}}
\newcommand\Supp{\operatorname{Supp}}
\newcommand\WCl{\operatorname{WCl}}

 \title{Embeddings of curves and surfaces
 \thanks{Research carried out under the EU HCM project AGE (Algebraic Geometry
in Europe), contract number ERBCHRXCT 940557. The final version was written
while the first author was ``{\em Professore distaccato}'' at the {\em
Accademia dei Lincei}.}}

 \author{F. Catanese \and M. Franciosi\thanks{Requests for preprints should
be addressed to the second author.} \and K. Hulek \and M. Reid}
 \date{Wed 5th Jun 1996}
 
\begin{document}
 \maketitle

 \begin{abstract} We prove a general embedding theorem for Cohen--Macaulay
curves (possibly nonreduced), and deduce a cheap proof of the standard results
on pluricanonical embeddings of surfaces, assuming vanishing $H^1(2K_X)=0$.
 \end{abstract}

 \section{Introduction}
Let $C$ be a curve over a field $k$ of characteristic $p\ge0$, and $H$ a
Cartier divisor on $C$. We assume that $C$ is projective and Cohen--Macaulay
(but possibly reducible or nonreduced). Write $HC=\deg\Oh_C(H)$ for the degree
of $H$, $p_aC=1-\chi(\Oh_C)$ for the arithmetic genus of $C$, and $\om_C$ for
the dualising sheaf (see \cite{Ha}, Chap.~III, \S7).

Our first result is the following. (A {\em cluster} $Z$ of {\em degree}
$\deg Z=r$ is simply a \hbox{$0$-dimensional} subscheme with
$\length\Oh_Z=\dim_k\Oh_Z=r$; a curve $B$ is {\em generically Gorenstein} if,
outside a finite set, $\om_B$ is locally isomorphic to $\Oh_B$. The remaining
definitions and notation are explained below.)

 \begin{TEO}[Curve embedding theorem]\label{th:curve} $H$ is very ample on $C$
if for every generically Gorenstein subcurve $B\subset C$, either
 \begin{enumerate}
 \item $HB\ge2p_aB+1$, or
 \item $HB\ge2p_aB$, and there does not exist a cluster $Z\subset B$ of
degree $2$ such that $\sI_Z\Oh_B(H)\iso\om_B$.
 \end{enumerate}

More generally, suppose that $Z\subset C$ is a cluster (of any degree) such
that the restriction
 \begin{equation}
H^{0}(C,\Oh_C(H))\to\Oh_Z(H)=\Oh_C(H)\otimes\Oh_Z
 \label{eq:rest}
 \end{equation}
is not onto. Then there exists a generically Gorenstein subcurve $B$ of $C$
and an inclusion $\fie\colon\sI_Z\Oh_B(H)\into\om_B$ not induced by a map
$\Oh_B(H)\to\om_B$. In particular, (\ref{eq:rest}) is onto if
 \begin{equation}
HB>2p_aB-2+\deg(Z\cap B)
 \nonumber
 \end{equation}
for every generically Gorenstein subcurve $B\subset C$.
 \end{TEO}

Theorem~\ref{th:curve} is well known for nonsingular curves $C$. Although
particular cases were proved in \cite{Ca1}, \cite{Ba2}, \cite{C-F}, \cite{C-H},
it was clear that the result was more general. In discussion after a lecture
on the Gorenstein case by the first author at the May 1994 Lisboa AGE meeting,
the fourth author pointed out the above result, where $C$ is only assumed to
be a pure 1-dimensional scheme. For divisors on a nonsingular surface, Mendes
Lopes \cite{ML} has obtained results analogous to Theorem~\ref{th:curve} and to
Theorem~\ref{th:hh}. We apply these ideas to the canonical map of a Gorenstein
curve in \S\ref{sec:cc}.

The proof of Theorem~\ref{th:curve} is based on two ideas from Serre and
Grothendieck duality:
 \begin{enumerate}
 \renewcommand\labelenumi{(\alph{enumi})}
 \item we use Serre duality in its ``raw'' form
 \begin{equation}
H^1(C,\sF)\dual\Hom(\sF,\om_C)\quad\text{for $\sF$ a coherent sheaf,}
 \nonumber
 \end{equation}
where $\dual$ denotes duality of vector spaces.
 \item If $\Oh_C$ has nilpotents, a nonzero map $\fie\colon\sF\to\om_C$ is
not necessarily generically onto; however (because we are $\Hom$'ming into
$\om_C$), duality gives an automatic factorisation of $\fie$ of the form
 \begin{equation}
\sF\onto\sF_{|B}\to\om_B\into\om_C,
 \nonumber
 \end{equation}
via a purely 1-dimensional subscheme $B\subset C$, where $\sF_{|B}\to\om_B$ is
generically onto. See Lemma~\ref{lem:adj} for details.
 \end{enumerate}

Since our main result might otherwise seem somewhat abstract and useless, we
motivate it by giving a short proof in \S\ref{sec:pluri}, following the methods
of \cite{C-F}, of the following result essentially due to Bombieri (when
$\chara k=0$) and to Ekedahl and Shepherd-Barron in general. Recall that a
{\em canonical surface} (or canonical model of a surface of general type) is a
surface with at worst Du Val singularities and $K_X$ ample. The remaining
notation and definitions are explained below.

 \begin{TEO}[Canonical embeddings of surfaces]\label{th:surf} $X$ is a
canonical surface. Assume that $H^1(2K_X)=0$. Then $mK_X$ is very ample if
$m\ge5$, or if $m=4$ and $K_X^2\ge2$, or if $m=3$, $p_g\ge2$ and $K_X^2\ge3$.
 \end{TEO}

Here $H^1(2K_X)=0$ follows at once in characteristic 0 from Kodaira vanishing
or Mumford's vanishing theorem. One can also get around the assumption
$H^1(2K_X)=0$ in characteristic $p>0$ (see \cite{Ek} or \cite{S-B}). In fact
Ekedahl's analysis (see \cite{Ek}, Theorem~II.1.7) shows that $H^1(2K_S)\ne0$
is only possible in a very special case, when $p=2,\chi(\Oh_S)=1$ and $S$ is
(birationally) an inseparable double cover of a K3 surface or a rational
surface.

In \S\ref{sec:tri} and \S\ref{sec:bi} we apply these ideas to prove the following theorems on
tricanonical and bicanonical linear systems of a surface of general type.

 \begin{TEO}[Tricanonical embeddings]\label{th:tri} Suppose that $X$ is a
canonical surface with $K_X^2\ge3$. Then $3K_X$ is very ample if either
 \begin{enumerate}
 \renewcommand\labelenumi{(\alph{enumi})}
 \item $q=h^1(\Oh_X)=0$; or
 \item $\chi(\Oh_X)\ge1$, $\dim\Pico X>0$ and $H^1(2K_X-L)=0$ for all
$L\in\Pico X$.
 \end{enumerate}
Note that (a) or (b) cover all cases with $\chara k=0$. Thus the cases {\em
not} covered by our argument are in $\chara k=p>0$, with either $p_g<q$ or
$\dim\Pico X=0$. \end{TEO}

Theorem~\ref{th:tri} in characteristic 0 is a result of Reider \cite{Rei}, but
see also \cite{Ca2}. Without the condition $K_X^2\ge3$, the double plane with
branch curve of degree 8 (that is, $X_8\subset\proj(1,1,1,4)$) is a
counter\-example. It follows from a result of Ekedahl (\cite{Ek},
Theorem~II.1.7) that if $\chi(\Oh_X)\ge1$ then $H^1(2K_X-L)=0$ for all $L\ne0$.
The remaining assumption in Theorem~\ref{th:tri} is that $H^1(2K_X)=0$, and
this can also be got around, as shown by Shepherd-Barron \cite{S-B}.

 \begin{TEO}[Bicanonical embeddings]\label{th:bi} We now assume that $q=0$
and $p_g\ge4$.

 \begin{enumerate}
 \renewcommand\labelenumi{(\alph{enumi})}
 \item $2K_X$ is very ample if every $C\in|K_X|$ is numerically $3$-connected
(in the sense of Definition~\ref{def:m-conn}, see also Lemma~\ref{lem:n-conn}).
More precisely, $|2K_X|$ separates a cluster $Z$ of degree $2$ provided that
every curve $C\in|K_X|$ through $Z$ is $3$-connected.
 \item Assume in addition that $K_X^2\ge10$, and let $Z$ be a cluster of
degree $2$ contained in $X$. Then $Z$ is contracted by $|2K_X|$ if and only if
$Z$ is contained in a curve $B\subset X$ with
 \begin{equation}
 K_XB=p_aB=1\text{ or }2
 \nonumber
 \end{equation}
(a {\em Francia curve}, compare Definition~\ref{def:Frc}), and
$\sI_Z\Oh_B(2K_X)\iso\om_B$.
 \item In particular, $|2K_X|$ defines a birational morphism unless $X$ has a
pencil of curves of genus $2$.
 \end{enumerate}
 \end{TEO}

 \begin{REMS}
(1) A cluster $Z$ of degree 2 is automatically contracted by $|2K_X|$ if it
is contained in a curve $C\subset X$ for which $\sI_Z\Oh_C(2K_X)\iso\om_C$ (for
a non\-singular curve, this reads ${2K_X}_{|C}=K_C+P+Q$). Thus (b) says in
particular that if this happens for some $C$ then it also happens for a
Francia curve.

(2) The assumptions $q=0$ and $p_g\ge4$ are needed for the simple minded
``restriction method'' of this paper, but we conjecture that (b) holds without
them (at least in characteristic zero, or assuming $q=0$); the case
$Z=\{x,y\}$ with $x\ne y$ (that is, ``separating points'') follows in
characteristic zero by Reider's method. We believe that the conjecture can be
proved quite generally by a different argument based on Ramanujam--Francia
vanishing, or by Reider's method applied to reflexive sheaves on $X$. Stay
tuned!

(3) In characteristic 0, Theorem~\ref{th:bi} (without the assumption $q=0$) is
due essentially to Francia (unpublished, but see \cite{Fr1}--\cite{Fr2}) and
Reider \cite{Rei}. Theorem~\ref{th:bi}, (a) is a consequence of
Theorem~\ref{th:hh} on canonical embeddings of curves and the generalisation
of hyper\-elliptic curves. The results in Theorem~\ref{th:bi} are only a
modest novelty, in that there is no restriction on the characteristic of the
ground field (see \cite{S-B}, Theorems~25, 26 and~27 for $\chara k\ge11$).

Further results on the bicanonical map $\fie_{2K}$ for smaller values of
$p_g$, $K_X^2$ (in characteristic 0) require a more intricate analysis, and we
refer to recent or forthcoming articles (\cite{C-F-M}, \cite{C-C-M}). Other
applications of our methods can be found in \cite{F}.
 \end{REMS}

 \subsection*{Acknowledgment} It is a pleasure to thank Ingrid Bauer for
interesting discussions on linear systems on surfaces, out of which this paper
originated.

 \subsection*{Conventions}
This paper deals systematically with reducible and nonreduced curves and their
subschemes $B\subset C$. A coherent sheaf $\sF$ on a curve $C$ is {\em torsion
free} if there are no sections $s\in\sF$ supported at points; on a
1-dimensional scheme, this is obviously equivalent to $\sF$ {\em
Cohen--Macaulay}. We say that $C$ is {\em purely \hbox{$1$-dimensional}\/} or
{\em Cohen--Macaulay} if $\Oh_C$ is torsion free.

A map $\fie\colon\sF\to\sG$ between coherent sheaves on $B$ is {\em
generically injective} if it is injective at every generic point of $B$; if
$\sF$ is torsion free then $\fie$ is automatically an inclusion $\sF\into\sG$.
If we know that the generic stalks of $\sF$ and $\sG$ have the same length at
every generic point of $C$ then a generically injective map
$\fie\colon\sF\to\sG$ is an isomorphism at each generic point, and therefore
$\ker\fie$ and $\coker\fie$ have finite length. Indeed, they are both coherent
sheaves supported at a finite set, and by the Nullstellensatz, each stalk is
killed by a power of the maximal ideal. This applies, for example, to the map
$\fie\colon\sI_Z\Oh_B(H)\into\om_B$ of Theorem~\ref{th:curve}, see
Lemma~\ref{lem:gg} below.

A scheme $B$ is {\em Gorenstein in codimension $0$} or {\em generically
Gorenstein} if $\om_B$ is locally isomorphic to $\Oh_B$ at every generic point
of $B$.

A {\em cluster} of degree $r$ is a 0-dimensional subscheme $Z\subset X$
supported at finitely many points, with ideal sheaf $\sI_Z$, structure sheaf
$\Oh_Z=\Oh_X/\sI_Z$, and having $\deg Z=h^0(\Oh_Z)=\length\Oh_Z=r$.
We sometimes write $Z=(x,y)$ for a cluster of degree 2, where $x,y$ are
either 2 distinct points of $X$, or a point $x$ plus a tangent vector $y$ at
$x$. We say that a linear system $|H|$ on $X$ {\em separates} $Z$ (or
separates $x$ and $y$) if $H^0(X,\Oh_X(H))\to\Oh_Z(H)$ is onto, or {\em
contracts} $Z$ if $Z$ does not meet the base locus $\Bs|H|$, and
$\rank\{H^0(X,\Oh_X(H))\to\Oh_Z(H)\}=1$.

 \subsection*{Notation}

 \begin{enumerate}

 \item[$X$] A projective scheme over an arbitrary field $k$. We sometimes (not
always consistently) write $k\subset\overline k$ for the algebraic closure,
and $X_{\overline k}=X\otimes_k\overline k$.

 \item[$\om_X$] Dualising sheaf of $X$ (see \cite{Ha}, Chap.~III, \S7).

 \item[$|H|$] Linear system defined by a Cartier divisor $H$ on $X$.

 \item[$C$] A curve, that is, a projective scheme over $k$ which is purely
1-dimensional, in the sense that $\Oh_C$ is Cohen--Macaulay (torsion free).

 \item[$p_aC$] The arithmetic genus of $C$, $p_aC=1-\chi(\Oh_C)$.

 \item[$K_C$] A canonical divisor of a Gorenstein curve $C$, that is, a Cartier
divisor such that $\Oh_C(K_C)\iso\om_C$ (only defined if $C$ is Gorenstein).

 \item[$\deg\sL$] The degree of a torsion free sheaf of rank 1 on $C$; it can be
defined by
 \begin{equation}
 \deg\sL=\chi(\sL)-\chi(\Oh_C).
 \nonumber
 \end{equation}
If $H$ is a Cartier divisor on $C$, we set $HC=\deg\Oh_C(H)$.

 \item[$S$] A nonsingular projective surface.

 \item[$DD'$] Intersection number of divisors $D,D'$ on a nonsingular
projective surface.

 \item[$K_S$] A canonical divisor on $S$.

 \item[$K_X^2$] If $X$ is a Gorenstein surface, $K_X$ is a Cartier divisor with
$\om_X=\Oh_X(K_X)$, and $K_X^2$ is the selfintersection number of the Cartier
divisor $K_X$. If $X$ has only Du Val singularities and $\pi\colon S\to X$
is the minimal nonsingular model then $K_S=\pi^*K_X$ and $K_X^2=K_S^2$.

 \item[$p_g,q$] The geometric genus $p_g=h^0(S,K_S)=h^0(X,K_X)$ of $S$ or $X$
(respectively the irregularity $q=h^1(S,\Oh_S)=h^1(X,\Oh_X)$).

 \item[$P_n$] The $n$th plurigenus $P_n=h^0(S,nK_S)$ of $S$.

 \end{enumerate}

 \section{Embedding curves} We start with a useful remark.

 \begin{REM}\label{rem:1} Let $H$ be a Cartier divisor on a scheme $X$. Then
$H$ is very ample if and only if the restriction map
 \begin{equation}
H^0(\Oh_X(H))\to\Oh_Z(H)
 \label{eq1}
 \end{equation}
is onto for every cluster $Z\subset X$ (more precisely,
for every $Z\subset X_{\overline k}$) of degree $\le2$.
 \end{REM}

 \begin{pf} By the standard embedding criterion of \cite{Ha}, Chap.~II,
Prop.~7.3, we have to prove that (\ref{eq1}) is onto for all the ideals
$\sI_Z=m_x$ or $m_xm_y$ with $x,y\in X$. For $x\ne y$, we are done.

By assumption $H^0(\Oh_C(H))\to\Oh_C/m_x$ is onto for every $x\in X$. Now if
the image of $H^0(m_x\Oh_C(H))\to m_x/m_x^2$ is contained in a hyperplane
$V\subset m_x/m_x^2$, then the inverse image of $V$ in $\Oh_{C,x}$ generates an
ideal $\sI\subset\Oh_{X,x}$ defining a cluster $Z$ of degree 2 supported at
$x$ such that $H^0(\Oh_C(H))\to\Oh_Z$ is not onto. \QED \end{pf}

 \begin{REM}\label{rem:2} The chain of reasoning we use below is that, by
Remark~\ref{rem:1} and cohomology, $H$ is very ample if and only if
$H^1(\sI_Z\Oh_X(H))\to H^1(\Oh_X(H))$ is injective for each cluster $Z$ of
degree $2$, or dually (if $X=C$ is a curve),
$\Hom(\Oh_C(H),\om_C)\to\Hom(\sI_Z(H),\om_C)$ is onto.
 \end{REM}

 \begin{LEM}\label{lem:gg} Let $C$ be a curve. Assume that there is a Cartier
divisor $H$ on $C$ and a cluster $Z\subset C$ for which the sheaf
$\sL=\sI_Z\Oh_C(H)$ has an inclusion $\sL\into\om_C$. Then $C$ is
generically Gorenstein.
 \end{LEM}

 \begin{pf} By assumption, $\sL\iso\Oh_C$ at every generic point of $C$. We must
prove that an inclusion $\sL\into\om_C$ maps onto every generic stalk
$\om_{C,\eta}$, or equivalently, that the cokernel $\sN=\om_C/\sL$ has finite
length. We give two slightly different proofs, one based on RR, and one using
properties of dualising modules.

Let $\Oh_C(1)$ be an ample line bundle on $C$. Then by Serre vanishing (see
\cite{Se1}, n$^\circ$~66, Theorem~2 or \cite{Ha}, Chap.~III, Theorem~5.2), for
$n\gg0$, the exact sequence
 \begin{equation}
0\to\sL(n)\to\om_C(n)\to\sN(n)\to0
 \nonumber
 \end{equation}
is exact on global sections, and all the $H^1$ vanish. Now by RR and duality,
 \begin{equation}
h^0(\om_C(n))=h^1(\Oh_C(-n))=-\chi(\Oh_C)+n\deg\Oh_C(1)\quad\text{for $n\gg0$.}
 \nonumber
 \end{equation}
On the other hand, RR also gives
$h^0(\sL(n))=\chi(\Oh_C)+\deg\sL+n\deg\Oh_C(1)$ for $n\gg0$, since
$\sL\iso\Oh_C$ at every generic point. Thus
 \begin{equation}
h^0(\sN(n))=-2\chi(\Oh_C)+\deg\sL\quad\text{for all $n\gg0$,}
 \nonumber
 \end{equation}
and therefore $\sN$ has finite length.

The alternative proof of the lemma uses the ``well-known fact'' (see below)
that the generic stalk $\om_{C,\eta}$ of the dualising sheaf at a generic
point $\eta\in C$ is a dualising module for the Artinian local ring
$\Oh_{C,\eta}$, so that they have the same length, and therefore an inclusion
$\sL\into\om_C$ is generically an isomorphism. The above proof in effect
deduces $\length\om_{C,\eta}=\length\Oh_{C,\eta}$ from RR together with Serre
duality, the defining property of $\om_C$.

 \paragraph{Proof of the ``well-known fact''} This is proof {\em by
incomprehensible reference}. First, if $\eta\in X$ is a generic point of a
scheme, more-or-less by definition, a dualising module of the Artinian ring
$\Oh_{X,\eta}$ is an injective hull of the residue field
$\Oh_{X,\eta}/m_{X,\eta}=k(\eta)$ (see \cite{Gr-Ha}, Proposition~4.10); in
simple-minded terms, $\Oh_{X,\eta}$ clearly contains a field $K_0$ such that
$K_0\subset k(\eta)$ is a finite field extension, and the vector space dual
$\Hom_{K_0}(\Oh_{X,\eta},K_0)$ is a dualising module. Next, if $\eta\in X$ is
a generic point of a subscheme $X\subset\proj=\proj^N$ of pure codimension $r$,
then by \cite{Ha}, Chap.~III, Prop.~7.5, the dualising sheaf of $X$ is
$\om_X=\sExt^r_{\Oh_\proj}(\Oh_X,\om_\proj)$. On the other hand, the local
ring $\Oh_{\proj,\eta}$ of projective space along $\eta$ is an
\hbox{$r$-dimensional} regular local ring, and therefore Gorenstein, so that
by \cite{Gr-Ha}, Prop.~4.13,
$\Ext^r_{\Oh_{\proj,\eta}}(\Oh_{X,\eta},\om_{\proj,\eta})$ is a dualising
module of $\Oh_{X,\eta}$ (an injective hull of the residue field
$\Oh_{X,\eta}/m_{X,\eta}=k(\eta)$). \QED \end{pf}

 \begin{LEM}[Automatic adjunction]\label{lem:adj} Let $\sF$ be a coherent
sheaf on $C$, and $\fie\colon\sF\to\om_C$ a map of $\Oh_C$-modules. Set
$\sJ=\Ann\fie\subset\Oh_C$, and write $B\subset C$ for the subscheme defined by
$\sJ$. Then $\fie$ has a canonical factorisation of the form
 \begin{equation}
\sF\onto\sF_{|B}\to\om_B=\sHom_{\Oh_C}(\Oh_B,\om_C)\subset\om_C,
 \label{eq:adj}
 \end{equation}
where $\sF_{|B}\to\om_B$ is generically onto.
 \end{LEM}

 \begin{pf} By construction of $\sJ$, the image of $\fie$ is contained in the
submodule
 \begin{equation}
\bigl\{s\in\om_C\bigm|\sJ s=0\bigr\}\subset\om_C
 \nonumber
 \end{equation}
But this clearly coincides with $\sHom(\Oh_B,\om_C)$. Now the inclusion
morphism $B\into C$ is finite, and $\om_B=\sHom_{\Oh_C}(\Oh_B,\om_C)$ is just
the adjunction formula for a finite morphism (see, for example, \cite{Ha},
Chap.~III, \S7, Ex.~7.2, or \cite{Re}, Prop.~2.11).

The factorisation (\ref{eq:adj}) goes like this: $\fie$ is killed by $\sJ$, so
it factors via the quotient module $\sF/\sJ\sF=\sF_{|B}$. As just observed, it
maps into $\om_B\subset\om_C$. Finally, it maps onto every generic stalk of
$\om_B$, again by definition of $\sJ$: a submodule of the sum of generic stalks
$\bigoplus\om_{B,\eta}$ is the dual to the generic stalk
$\bigoplus\Oh_{B',\eta}$ of a purely 1-dimensional subscheme $B'\subset B$,
and $\fie$ is not killed by the corresponding ideal sheaf $\sJ'$. \QED
 \end{pf}

 \begin{REM} We define $B$ to be the {\em scheme theoretic support} of $\fie$.
Note that if $C=\sum n_i\Ga_i$ is a Weil divisor on a normal surface and $\sF$
a line bundle, the curve $B\subset C$ defines a splitting $C=A+B$ where $A$ is
the {\em divisor of zeros} of $\fie$: at the generic point of $\Ga_i$, the
map $\fie$ then looks like $y_i^{a_i}$, where $y_i$ is the local equation
of $\Ga_i$, and $A=\sum a_i\Ga_i$. In the general case however, $A$ does not
make sense. \end{REM}

 \begin{pfof}{Theorem~\ref{th:curve}} Let $H$ be a Cartier divisor, and $\sI$
the ideal sheaf of a cluster for which $H^1(\sI\Oh_C(H))\ne0$. Then
$\Hom(\sI\Oh_C(H),\om_C)\ne0$ by Serre duality. First pick any nonzero map
$\fie\colon\sI\Oh_C(H)\to\om_C$. By Lemma~\ref{lem:adj}, $\fie$ comes from an
inclusion $\sI\Oh_B(H)\into\om_B$ for a subscheme $B\subset C$, and $B$ is
generically Gorenstein by Lemma~\ref{lem:gg}.

Finally, if $H^0(\Oh_C(H))\to\Oh_Z(H)$ is not onto, then the next arrow in
the cohomology sequence
 \begin{equation}
H^1(\sI\Oh_C(H))\to H^1(\Oh_C(H))
 \nonumber
 \end{equation}
is not injective, and dually, the restriction map
 \begin{equation}
\Hom(\Oh_C(H),\om_C)\to\Hom(\sI\Oh_C(H),\om_C)
 \nonumber
 \end{equation}
is not onto. Thus we can pick $\fie\colon\sI\Oh_C(H)\to\om_C$ which is not the
restriction of a map $\Oh_C(H)\to\om_C$. Then also the map
$\sI\Oh_B(H)\into\om_B$ given by Lemma~\ref{lem:adj} is not the restriction of
a map $\Oh_B(H)\into\om_B$.

For the final part, an inclusion $\sI\Oh_B(H)\into\om_B$ has cokernel of
finite length, so that $\chi(\sI\Oh_B(H))\le\chi(\om_B)$. Plugging in the
definition of degree gives
 \begin{equation}
1-p_aB+HB-\deg(Z\cap B)\le p_aB-1,
 \nonumber
 \end{equation}
that is,
 \begin{equation}
HB\le2p_aB-2+\deg(Z\cap B).
 \nonumber
 \end{equation}
Thus, assuming the inequality (2) of Theorem~\ref{th:curve}, no such inclusion
$\sI\Oh_B(H)\into\om_B$ can exist, so that $H^0(\Oh_C(H))\to\Oh_Z(H)$ is
onto. \QED \end{pfof}

 \section{The canonical map of a Gorenstein curve}\label{sec:cc} We now discuss
the canonical map $\fie_{K_C}$ of a Gorenstein curve, writing $K_C$ for a
canonical divisor of $C$, that is, a Cartier divisor for which
$\om_C\iso\Oh_C(K_C)$. Our approach is motivated in part by the examples and
results in the reduced case treated in \cite{Ca1}.

 \begin{DEF}\label{def:m-conn} A Gorenstein curve $C$ over an algebraically
closed field
$k$ is {\em numerically $m$-connected} if
 \begin{equation}
\deg\Oh_B(K_C)-\deg\om_B=\deg(\om_C\otimes\Oh_B)-(2p_aB-2)\ge m
 \nonumber
 \end{equation}
for every generically Gorenstein strict subcurve $B\subset C$. For $C$ over
any field, we say that $C$ is numerically $m$-connected if
$C\otimes\overline k$ is numerically $m$-connected.
 \end{DEF}

 \begin{REM}\label{rem:m-conn} Note that for divisors on a nonsingular
surface,
 \begin{equation}
\deg\Oh_B(K_C)-\deg\om_B=(K_S+C)B-(K_S+B)B=(C-B)B.
 \nonumber
 \end{equation}
In this context, Franchetta and Ramanujam define numerically connected in
terms of the intersection numbers $AB=(C-B)B$ for all effective decompositions
$C=A+B$. The point of our definition is to use the numbers
$\deg\Oh_B(K_C)-\deg\om_B$ in the more general case as a substitute for
$(C-B)B$. In effect, we think of the adjunction formula as defining the
``degree'' of the ``normal bundle'' to $B$ in $C$, in terms of the difference
between $K_C{}_{|B}$ and $\om_B$. \end{REM}

 \begin{TEO}\label{th:free} Let $C$ be a Gorenstein curve over a field $k$.
 \begin{enumerate}
 \renewcommand\labelenumi{(\alph{enumi})}
 \item If $C$ is numerically $1$-connected then
$H^0(\Oh_C)=k$ (the constant functions).

 \item If $C$ is numerically $2$-connected then either $|K_C|$ is free or
$C\iso\proj^1$ (over the algebraic closure $\overline k$, of course). In
particular, in this case $p_aC=0$ implies $C\iso\proj^1$.
 \end{enumerate}
 \end{TEO}
 \begin{pfof}{(a)} First, if $f\in H^0(\Oh_C)$ is a nonzero section vanishing
along some reduced component of $C$, then applying Lemma~\ref{lem:adj} to the
multiplication map $\mu_f\colon\Oh_C(K_C)\to\om_C$ gives an inclusion
$\Oh_B(K_C)\into\om_B$, which is forbidden by numerically 1-connected (because
$\deg\Oh_B(K_C)>\deg\om_B$). Now if $H^0(\Oh_C)\ne k$, there exists a nonzero
section $f\in H^0(\Oh_{C\otimes\overline k})$ vanishing at any given point
$x\in C\otimes\overline k$. An inclusion $\Oh_C\into m_x$ contradicts at once
$0=\deg\Oh_C>\deg m_x=-1$, so that $f$ must vanish along some component of
$C$, and we have seen that this is impossible. \QED \end{pfof}

 \begin{pfof}{(b)} As discussed in Remark~\ref{rem:2}, the standard chain of
reasoning is as follows:
 \begin{enumerate}
 \item $x\in C$ is a base point of $|K_C|$ if and only if
$H^0(\Oh_C(K_C))\to\Oh_x(K_C)$ is not onto, and then
 \item $H^1(m_x\Oh_C(K_C))\to H^1(\Oh_C(K_C))$ is not injective,
 \item dually, $\Hom(\Oh_C(K_C),\om_C)\to\Hom(m_x\Oh_C(K_C),\om_C)$ is not
onto,
 \item therefore there exists a map $s\colon m_x\Oh_C(K_C)\to\om_C$ linearly
independent of the identity inclusion.
 \end{enumerate}

Now by Lemma~\ref{lem:adj}, the map $s$ factors via an inclusion
$m_x\Oh_B(K_C)\into\om_B$ on a generically Gorenstein curve $B$. But then
$B\subsetneq C$ is forbidden by the numerically 2-connected assumption
$\deg m_x\Oh_B(K_C)-\deg\om_B\ge1$.

Therefore $B=C$, that is, $s\colon m_x\Oh_C(K_C)\into\om_C$ is an inclusion.
After tensoring down with $-K_C$, this gives an inclusion
$i\colon m_x\into\Oh_C$ linearly independent of the identity. Write
$\sF=i(m_x)\subset\Oh_C$. Then $\deg\sF=-1$, and therefore $\sF=m_z$ for some
$z\in C$.

Now for any point $y\in C\setminus\{x\}$, there exists a linear combination
$s'=s+\la\id$ vanishing at $y$, which therefore defines an isomorphism
$m_x\iso m_y$. This implies that every point $y\in C$ is a Cartier divisor,
hence a nonsingular point. Since $C$ is clearly connected, and
$\Oh_C(x-y)\iso\Oh_C$ for every $x,y\in C$, it follows that $C\iso\proj^1$.

For the final statement, if $p_aC=0$ then $1=h^0(\Oh_C)=h^1(\om_C)$ by
(a) and duality, hence $h^0(\om_C)=0$ by RR, so that $H^0(\Oh_C(K_C))\to\Oh_x$
is not onto for any $x\in C$. \QED \end{pfof}

 \begin{DEF} We say that a Gorenstein curve $C$ is {\em honestly hyperelliptic
(\cite{Ca1}, Definition~3.18)} if there exists a finite morphism $\psi\colon
C\to\proj^1$ of degree $2$ (that is, $\psi$ is finite and $\psi_*\Oh_C$ is
locally free of rank $2$ on $\proj^1$). The linear system $\psi^*|\Oh_C(1)|$
defining $\psi$ is called an {\em honest $g^1_2$.} \end{DEF}

We note the immediate consequences of the definition.

 \begin{LEM}\label{lem:hh} An honestly hyperelliptic curve $C$ of genus
$p_aC=g\ge0$ is isomorphic to a divisor $C_{2g+2}$ in the weighted projective
space $\proj(1,1,g+1)$, not passing through the vertex $(0,0,1)$, defined by an
equation
 \begin{equation}
w^2+a_{g+1}(x_1,x_2)w+b_{2g+2}(x_1,x_2)=0.
 \nonumber
 \end{equation}
It follows that every point of $C$ is either nonsingular or a planar double
point, and that $C$ is either irreducible, or of the form $C=D_1+D_2$ with
$D_1D_2=g+1$.

The projection $\fie\colon C\to\proj^1$ is a finite double cover, and the
inverse image of any $x\in\proj^1$ is a Cartier divisor which is a cluster
$Z\subset C$ of degree $2$. In other words, $Z$ is either $2$ distinct
nonsingular points of $C$, a nonsingular point with multiplicity $2$, or a
section through a planar double point of $C$. \qed\end{LEM}

 \begin{TEO}\label{th:hh} Let $C$ be a numerically $3$-connected Gorenstein
curve. Then either $|K_C|$ is very ample or $C$ is honestly hyperelliptic.

In particular, in this case if\/ $p_aC\ge2$ then $K_C$ is ample, and if\/
$p_aC=1$ then $C$ is honestly hyperelliptic (over the algebraic closure
$\overline k$, of course).
 \end{TEO}

 \begin{pf} Let $Z$ be a cluster of degree 2 for which
$H^0(\Oh_C(K_C))\to\Oh_Z(K_C)$ is not onto. The previous chain of reasoning
gives a map $\sI_Z\Oh_C(K_C)\to\om_C$ linearly independent of the identity
inclusion. An inclusion $\sI_Z\Oh_B(K_C)\into\om_B$ with $B\subsetneq C$ is
forbidden as before by $C$ numerically 3-connected. Therefore we get an
inclusion $s\colon\sI_Z\Oh_C(K_C)\into\om_C$ linearly independent of the
identity inclusion. Note that any linear combination $s'=s+\la\id$ of the two
sections is again generically injective, since an inclusion
$\sI_Z\Oh_B(K_C)\into\om_B$ with $B\subsetneq C$ is forbidden by numerically
3-connected.

The image $\sF=s(\sI_Z\Oh_C(K_C))\subset\om_C$ is a submodule of colength 2,
therefore of the form $\sF=\sI_{Z'}\om_C$ for some cluster $Z'\subset C$.
Tensoring down the iso\-morphism $s\colon\sI_Z\Oh_C(K_C)\to\sI_{Z'}\om_C$
gives an isomorphism $s\colon\sI_Z\iso\sI_{Z'}$, still linearly independent
of the identity inclusion $\sI_Z\into\Oh_C$.

Logically, there are 3 cases for $Z$ and $Z'$. The first of these
corresponds to an honest $g^1_2$ on $C$; the other two, corresponding to a
$g^1_2$ with one or two base points, lead either to $p_aC\le1$ or to a
contradiction. The case division is as follows:

 \paragraph{Case $Z\cap Z'=\emptyset$} Then the isomorphism
$\sI_Z\iso\sI_{Z'}$ implies that both $Z$ and $Z'$ are Cartier
divisors, and the two linearly independent inclusions $\sI_Z\into\Oh_C$
define an honest $g^1_2$ on $C$. In more detail: $\Oh_C(Z)$ has 2 linearly
independent sections with no common zeroes, and no linear combination of these
vanishes on any component of $C$. Therefore $|Z|$ defines a finite 2-to-1
morphism $C\to\proj^1$.

 \paragraph{Case $Z=Z'$} This case leads to an immediate contradiction.
Indeed, take any point $x\notin\Supp Z$; then some linear combination of the
two isomorphisms $s,\id\colon\sI_Z\to\sI_Z$ vanishes at $x$, and
therefore vanishes along any reduced component of $C$ containing $x$. But we
have just said that this is forbidden.

 \paragraph{Case $Z\cap Z'=x$} Here the case assumption can be rewritten
$\sI_Z+\sI_{Z'}=m_x$. This case is substantial, and it really happens in
two examples:
 \begin{enumerate}
 \item if $C$ is an irreducible plane cubic with a node or cusp $P$, and
$Q,Q'\in C\setminus P$ then $m_Pm_Q\iso m_Pm_{Q'}$;
 \item $\proj^1$ has an incomplete $g^1_2$ with a fixed point, of the form
$P+|Q|$.
 \end{enumerate}
We prove that we are in one of these cases. In either example, the curve $C$
has an honest $g^1_2$ (not directly given by our sections $s,\id$), so the
theorem is correct.

 \begin{CLA}\label{cla:mov_y} For any point $y\in C\setminus\{x\}$, there
exists a linear combination $s'=s+\la\id$ defining an isomorphism $\sI_Z\iso
m_xm_y$.
 \end{CLA}

 \begin{pfof}{Claim}
Since $\sI_Z,\sI_{Z'}\subset m_x$, we have two linearly independent maps
$s,\id\colon\sI_Z\into m_x$, and some linear combination $s'=s+\la\id$
vanishes at $y$. Also, no map $\sI_Z\to m_x$ vanishes along a component of
$C$. Thus $s'(\sI_Z)=m_xm_y$. \QED \end{pfof}

It follows from the claim that $m_xm_y\iso m_xm_{y'}$ for any two points
$y,y'\ne x$, so that $y,y'$ are nonsingular, and $C$ is reduced and
irreducible. Now let $\si\colon C_1\to C$ be the blow up of $m_x$. Then,
essentially by definition of the blow up, $m_x\Oh_{C_1}\iso\Oh_{C_1}(-E)$
where $E$ is a Cartier divisor on $C_1$. Then $m_{C_1,y}\iso m_{C_1,y'}$
for general points $y,y'\in C_1$, hence as usual $C_1\iso\proj^1$. If
$C_1\iso C$ there is nothing more to prove.

If $C_1\not\iso C$, the conductor ideal
$\sC=\sHom_{\Oh_C}(\si_*\Oh_{C_1},\Oh_C)$ of $\si_*\Oh_{C_1}$ in $\Oh_C$ is
$m_x$. Indeed, let $f\in k(C)$ be the rational function such that
multiplication by $f$ gives $m_xm_y\iso m_xm_{y'}$; then $f$ is an affine
parameter on $C_1=\proj^1$ outside $y$, so that all regular functions on $C_1$
are regular functions of $f$, and $fm_x=m_x$ implies
$\si_*(m_x\Oh_{C_1})=m_x\subset\Oh_C$. Now it is known that the only
Gorenstein curve singularity $x\in C$ with conductor ideal $m_x$ is a node or
cusp (see \cite{Se2}, Chap.~IV, \S11 or \cite{Re}, Theorem~3.2): indeed,
$m_x\subset\Oh_C\subset\si_*\Oh_{C_1}$, and the Gorenstein assumption $n=2\de$
gives $\length(\si_*\Oh_{C_1}/\Oh_C)=\length(\Oh_C/m_x)=1$. Therefore
$p_aC=1$.

For the final statement, if $p_aC=1$ then $1=h^0(\Oh_C)=h^1(\om_C)$ by
Theorem~\ref{th:free}, (a) and duality, hence $h^0(\om_C)=1$ by RR, so that
$H^0(\Oh_C(K_C))\to\Oh_Z$ is not onto for any cluster $Z\in C$ of degree
2. \QED \end{pf}

 \begin{REM} If $C$ is a numerically $3$-connected Gorenstein curve with
$p_aC\ge2$, then Theorem~\ref{th:hh} says that $K_C$ is automatically ample,
and the usual dichotomy holds: either $K_C$ is very ample, or $C$ is honestly
hyperelliptic.

Now assume instead that the dualising sheaf $\om_C=\Oh_C(K_C)$ is ample and
generated by its $H^0$. Equivalently, that $|K_C|$ is a free linear system,
defining a finite morphism (the {\em canonical morphism})
$\fie=\fie_{K_C}\colon C\to\proj^{p_a-1}$. In \cite{Ca1}, Definition~3.9, $C$
was defined to be {\em hyperelliptic} if $\fie_{K_C}$ is not birational on
some component of $C$. Thus by Theorem~\ref{th:hh}, in the $3$-connected case,
hyperelliptic and honestly hyperelliptic coincide.
 \end{REM}

 \section{Canonical maps of surfaces of general type}\label{sec:pluri} We give
a slight refinement of a useful lemma due independently to J.~Alexander and
I.~Bauer.

 \begin{LEM}[Alexander--Bauer]\label{lem:ab} Suppose that $H$ is a Cartier
divisor on an irreducible projective scheme $X$. Assume given effective Cartier
divisors $D_1,D_2$, $D_3$ such that

 \begin{enumerate}
 \renewcommand\labelenumi{(\roman{enumi})}
 \item $H^0(\Oh_X(H))\to H^0(\Oh_{D_i}(H))$ is onto.

 \item $H$ is very ample on every $\De\in|D_i|$ for $i=1,2,3$.
 \end{enumerate}

Then $H$ is very ample on $X$ if either
 \begin{enumerate}
 \renewcommand\labelenumi{(\alph{enumi})}
 \item $H\lineq D_1+D_2$ and $\dim|D_2|\ge1$, or

 \item $H\lineq D_1+D_2+D_3$ and $\dim|D_i|\ge1$ for $i=1,2,3$.
 \end{enumerate}
 \end{LEM}

 \begin{pf} (a) is proved in \cite{Ba1}, Claim~2.19 and \cite{Ra}, Lemma~3.1,
and also in \cite{C-F}, Prop.~5.1.

We prove (b). By Remark~\ref{rem:1}, we need to prove that if $x$ is any point
of $X$, and $y$ is either another point of $X$ or a tangent vector at $x$,
then $|H|$ separates $x$ from $y$. If some $\De_i\in|D_i|$ contains both $x$
and $y$, we are done by the assumptions (i) and (ii). In particular, since
$\dim|D_i|\ge1$, such a $\De_i$ exists if $x$ or $y$ belong to the base
locus of $|D_i|$.

Finally, if none of the above possibilities occurs, we can find $\De_1$
containing $x$ but not $y$, and $\De_2,\De_3$ containing neither $x$ nor $y$.
Then $\De_1+\De_2+\De_3$ separates $x$ from $y$. \QED\end{pf}

 \begin{pfof}{Theorem~\ref{th:surf}} Let $\pi\colon S\to X$ be the natural
birational morphism from a minimal surface of general type $S$ to its
canonical model $X$; write $K_S$ and $K_X$ for the canonical divisors of $S$
and $X$. Then $\om_X$ is invertible and $\pi^*(\om_X)\iso\om_S$; in particular
$H^0(X,mK_X)\iso H^0(X,mK_S)$ and$K_X^2=K_S^2$.

 \subparagraph{Step I} If $C\in|(m-2)K_X)|$, then $H^0(\Oh_X(mK_X))\to
H^0(\Oh_C(mK_X))$ is onto. This follows from our assumption
$H^1(\Oh_X(2K_X))=0$.

 \subparagraph{Step II} If $C\in|(m-2)K_X|$, then $\Oh_C(mK_X)$ is very ample.

 \begin{pf} By the curve embedding theorem Theorem~\ref{th:curve}, it is enough
to prove that $mK_XB\ge2p_aB+1$ for every subcurve $B\subset C$. Note that by
adjunction $K_C=(m-1){K_X}_{|C}$, so that we can write
$m{K_X}_{|B}={K_X}_{|B}+{K_C}_{|B}$. Since $K_X$ is ample,
$K_XB\ge1$, and therefore we need only prove that $K_XC\ge3$ and
 \begin{equation}
\deg\Oh_B(K_C)-\deg\om_B\ge2\quad\text{for every strict subcurve $B\subset C$,}
 \nonumber
 \end{equation}
that is, that $C$ is numerically 2-connected.

The corresponding fact for the minimal nonsingular model $S\to X$ is easy and
well known.\footnote{{\bf Tutorial}\enspace This is an easy consequence of the
Hodge algebraic index theorem. If $D$ is nef and big and $D=A+B$ then
$A^2+AB\ge0$, $AB+B^2\ge0$. The index theorem says that $A^2B^2\le(AB)^2$, with
equality only if $A,B$ are numerically equivalent to rational multiples of one
another. The reader should carry out the easy exercise of seeing that $AB\le0$
gives a contradiction, and proving all the connected assertions we need. Or
see \cite{Bo}, \S4, Lemma~2 for details (the exceptional case $n=2$,
$2K_S=A+B$, with $A\numeq B\numeq K_S$ and $K_S^2=1$ is
excluded by the assumption $K_S^2\ge2$ if $m=4$ of Theorem~\ref{th:surf}).}
Therefore $C$ numerically 2-connected follows from the next result, whose
proof we relegate to an appendix.

 \begin{LEM}\label{lem:n-conn} Let $X$ be a surface with only Du Val
singularities, and $\pi\colon S\to X$ the minimal resolution of singularities.
Let $C\subset X$ be an effective Cartier divisor, and $C^*=\pi^*C$ the total
transform of $C$ on $S$. Then
 \begin{equation}
\text{$C^*$ numerically $k$-connected} \implies \text{so is $C$.}
 \nonumber
 \end{equation}

Moreover, if $C^*$ is numerically $2$-connected, and is only $3$-disconnected
by expressions $C^*=A+B$ where $A$ or $B$ is a $-2$-cycle exceptional for
$\pi$ then $C$ is numerically $3$-connected.
 \end{LEM}\unskip\end{pf}

 \subparagraph{Step III} $h^0((m-2)K_X)\ge3$ if $m\ge5$, and $\ge2$ if $m=3$ or
4.

 \begin{pf} For $m=3$ this is just the assumption $p_g\ge2$.

For $m\ge4$, if $p_g\ge2$, then clearly $h^0((m-2)K_X)\ge3$. Otherwise, in
the case $p_g\le1$, we use the traditional numerical game of \cite{B-M},
based on Noether's formula $12\chi(\Oh_X)=(c_1^2+c_2)(X)$. It consists of
writing out Noether's formula using Betti numbers for the etale cohomology, in
the form
 \begin{equation}
 10+12p_g=8h^1(\Oh_X)+2\De+b_2+K_X^2.
 \label{eq:No}
 \end{equation}
Here the nonclassical term $\De=2h^1(\Oh_X)-b_1$ satisfies $\De\ge0$, and
$\De=0$ if $\chara k=0$. Since all the terms on the right hand side of
(\ref{eq:No}) are $\ge0$, it follows immediately that
 \begin{equation}
 \begin{aligned}
&p_g\le1\implies h^1(\Oh_X)\le2\\
&p_g\le0\implies h^1(\Oh_X)\le1.
 \end{aligned}
 \nonumber
 \end{equation}
Therefore, $p_g\le1$ implies $\chi(\Oh_X)\ge0$; hence, for $m\ge4$, by
RR
 \begin{equation}
h^0((m-2)K_X)\ge\chi(\Oh_X)+\binom{m-2}2K_X^2\quad
 \begin{cases}
\ge3&\text{if $m\ge5$,}\\
\ge2&\text{if $m=4$.}
 \end{cases}
 \nonumber
 \end{equation}
 \end{pf}
\unskip
 \subparagraph{Step IV} For $m=3$, we simply apply Lemma~\ref{lem:ab}, (b) to
$3K\lineq K+K+K$. For $m=4$ we apply Lemma~\ref{lem:ab}, (a) to $4K\lineq 2K+2K$:
the assumptions (i) and (ii) of the lemma hold by Steps~I, II and~III.

For $m\ge5$, we want to show that $H^0(\Oh_X(mK_X))\to\Oh_Z$ is onto for any
cluster $Z\subset X$ of degree 2. But by Step~III, there exists
$C\in|(m-2)K_X|$ containing $Z$. The result then follows by Steps~I and~II.
\QED \end{pfof}

 \subsection*{Appendix: Proof of Lemma~\ref{lem:n-conn}}

Suppose that $B\subset C$ is a strict subcurve. Write $B'$ for the birational
(=strict or proper) transform of $B$ in $S$ and $C^*=\pi^*C$ for the total
transform of $C$. For the proof, we find a divisor $\Bh$ (the {\em hat
transform}) with the properties
 \begin{enumerate}
 \renewcommand\labelenumi{(\roman{enumi})}
 \item $B'\le\Bh\le C^*$ and $\Bh-B'$ contains only exceptional curves;
 \item $p_a\Bh=p_aB$.
 \end{enumerate}
Suppose first that we know $\Bh$ satisfying these conditions. Then
 \begin{equation}
 (C^*-\Bh)\Bh\ge k
 \nonumber
 \end{equation}
by the assumption on $C^*$, which we write
 \begin{equation}
 (K_S+C^*)\Bh-(K_S+\Bh)\Bh\ge k.
 \nonumber
 \end{equation}
Here the first term equals $(K_X+C)B=\deg\Oh_B(K_C)$, and the second
$2p_a\Bh-2=2p_aB-2$. Thus
 \begin{equation}
\deg\Oh_B(K_C)-(2p_aB-2)=(K_S+C^*)\Bh-(2p_a\Bh-2)\ge k.
 \nonumber
 \end{equation}

So it is enough to find $\Bh$. For this, following the methods of
\cite{Ar1}--\cite{Ar2}, let $\bigl\{\Ga_i\bigr\}$ be all the exceptional
$-2$-curves. Define $\Bh=B'+\sum e_i\Ga_i$ with $e_i\in\Z$, $e_i\ge0$ minimal
with respect to the property $\Bh\Ga_i\le0$; this exists, because $C^*-A'$ has
the stated property (where $A'$ is the birational transform of the residual
Weil divisor $C-B$).

 \begin{CLA}\label{cla:bhat} The curve $\Bh$ has the following properties:
 \begin{enumerate}
 \renewcommand\labelenumi{(\roman{enumi})}
 \setcounter{enumi}2
 \item $\om_B=\pi_*\om_{\Bh}$;
 \item $R^1\pi_*\om_{\Bh}=0$.
 \end{enumerate}
Therefore $p_a\Bh=p_aB$.
 \end{CLA}

 \begin{pfof}{Claim} Taking $\pi_*$ of the short exact sequence
 \begin{equation}
0\to\Oh_S(K_S)\to\Oh_S(K_S+\Bh)\to\om_{\Bh}\to0
 \nonumber
 \end{equation}
gives $0\to\Oh_X(K_X)\to\Oh_X(K_X+B)\to\pi_*\om_{\Bh}\to0$ and
$R^1\pi_*\Oh_S(K_S+\Bh)=R^1\pi_*\om_{\Bh}$. The first of these implies that
$\om_B=\pi_*\om_{\Bh}$. Indeed, if $B\subset X$ is an effective Weil divisor
on any Cohen--Macaulay variety then the adjunction formula
$\om_B=\sExt^1_{\Oh_X}(\Oh_B,\om_X)$ (see, for example, \cite{Re},
Theorem~2.12, (1)) boils down to an exact sequence
$0\to\Oh_X(K_X)\to\Oh_X(K_X+B)\to\om_B\to0$. This proves (iii).

By the method of \cite{Ar1}--\cite{Ar2},
 \begin{equation}
R^1\pi_*\Oh_S(K_S+\Bh)=\varprojlim H^1(D,\Oh_D(K_S+\Bh)),
 \nonumber
 \end{equation}
where the inverse limit is taken over effective divisors $D=\sum a_j\Ga_j$. If
all the $H^1=0$, the limit is zero, as required.

Suppose then by contradiction that $D=\sum a_j\Ga_j$ has
$H^1(\Oh_D(K_S+\Bh))\ne0$. Then dually, $\Hom(\Oh_D(K_S+\Bh),\om_D)\ne0$, and
Lemma~\ref{lem:adj} gives an inclusion $\Oh_D(K_S+\Bh)\into\om_D$ (for a
possibly smaller $D$). Writing out the adjunction formula for $\om_D$ and
tensoring down by $K_S+\Bh$ gives $\Oh_D\into\Oh_D(D-\Bh)$. Therefore
$(\Bh-D)\Ga_i\le0$ for every $\Ga_i\subset D$, and by construction of $\Bh$
for the other $\Ga_i$. Now $\Bh-D=B'+\sum e'_j\Ga_j$ contradicts the
minimality of $\Bh$, provided we show that the $e_j'\ge0$. For this, note that
 \begin{equation}
\bigl(\sum e'_j\Ga_j\bigr)\Ga_i=(\Bh-D)\Ga_i-B'\Ga_i\le0
\quad\text{for every $i$}
 \nonumber
 \end{equation}
and the intersection form on the $\Ga_i$ is negative definite, so that the
standard argument implies $\sum e'_j\Ga_j\ge0$ (write it as $A-B$ where
$A,B\ge0$ have no common divisor, and calculate $B^2$). \QED \end{pfof}

 \section{The tricanonical map}\label{sec:tri}

We state the following three points as independent lemmas in order to tidy up
our proofs, and because they might be useful elsewhere. The first is a
particular case of the numerical criterion for flatness, see \cite{Ha},
Chap.~III, Theorem~9.9.

 \begin{LEM}[Flat double covers]\label{lem:flat} If $\fie\colon X\to Y$ is a
generically $2$-to-$1$ morphism (say with $Y$ integral), then for any $y\in Y$,
the condition $\length\fie\1(y)=2$ implies that $\fie$ is flat over a
neighbourhood of $y$. \qed \end{LEM}

 \begin{LEM}[Push-down of invariant linear systems]\label{lem:push} Let
$\fie\colon X\to Y$ be a finite morphism of degree $2$, where $X$ and $Y$ are
normal. Suppose that $L$ is a linear system of Cartier divisors on $X$ with
the property that $\fie_{|D}\colon D\to\Ga_D=\fie(D)$ has degree $2$ for every
$D\in L$. Then the $\Ga_D$ are linearly equivalent Weil divisors, that is, they
are all members of one linear system.
 \end{LEM}

 \begin{pf} For any $D,D'\in L$, note that $2\Ga_D=\pi_*D$ is a Cartier
divisor on $Y$, and $2\Ga_D\lineq 2\Ga_{D'}$, because if $D$ is locally
defined by $f\in k(X)$ (or $D-D'=\div f$) then $2\Ga_D$ is locally defined by
$\Norm(f)$, where $\Norm=\Norm_{k(X)/k(Y)}$.

Thus the Weil divisor class $\Ga_D-\Ga_{D'}$ is a 2-torsion element of the
Weil divisor class group $\WCl Y$ (modulo linear equivalence). The group of
Weil divisors numerically equivalent to zero is an algebraic group of finite
type, so that its 2-torsion subgroup is a finite algebraic group scheme $G$.
Now for fixed $D_0\in L$, taking $D\mapsto\Ga_D-\Ga_{D_0}$ defines a morphism
from the parameter space of the linear system $L$ to $G$, which must be the
constant morphism to 0. This proves what we need.

Assuming that $\fie$ is separable make this argument more intuitive, since
then it is Galois, and $\fie_*\Oh_X$ splits into invariant and antiinvariant
parts: $\fie_*\Oh_X=\Oh_Y\oplus\sL$, with $\sL$ a divisorial sheaf. Then
$\Ga_D$ is locally either a Cartier divisor or in the local Weil divisor class
of $\sL$, and $\Ga_D-\Ga_{D'}$ is in the kernel of $\fie^*$, which is a finite
algebraic group scheme, etc. \QED \end{pf}

 \begin{LEM}\label{lem:sing} Let $\La$ be a linear system of Weil divisors
through a point $P$ on a normal surface $Y$. Then the curves in $\La$ singular
at $P$ form a projective linear subspace of codimension $\le2$.
 \end{LEM}
 \begin{pf} Easy exercise involving the resolution and birational transform.
\qed \end{pf}

 \begin{pfof}{Theorem~\ref{th:tri}, Case~(a)} Since $q=0$, we have
$\chi(\Oh_X)\ge1$, and $K_X^2\ge3$ gives $P_2=h^0(2K_X)\ge4$. Let $Z$ be a
cluster of degree 2 on $X$. Since $P_2\ge4$, the linear subsystem
$|2K_X-Z|$ consisting of curves $D\in|2K_X|$ through $Z$ has dimension
$\ge1$, and any $D\in|2K_X|$ is 3-connected by the final part of
Lemma~\ref{lem:n-conn} (whose assumptions are easily verified as in \cite{Bo},
\S4, Lemma~2). By $H^1(K_X)=0$, the sequence
 \begin{equation}
0\to H^0(X,\Oh_X(K_X))\to H^0(X,\Oh_X(3K_X))\to H^0(D,\om_D)\to0
 \nonumber
 \end{equation}
is exact. Since $|\om_D|$ is free by Theorem~\ref{th:free}, it follows that
$\fie=\fie_{3K_X}$ is a finite morphism $\fie\colon X\to Y\subset\proj^N$,
where $N=P_3-1$. Assume that $|3K_X|$ does not separate $Z$. Then, by
Theorem~\ref{th:hh}, $D$ is honestly hyperelliptic. Since the same argument
applies to any $D\in|2K_X-Z|$, it follows that $\deg\fie\ge2$.

On the other hand, for any point $y\in Y$, if the scheme theoretic fibre
$\fie\1(y)$ is a cluster of degree $\ge3$, then there is a curve $D'\in|2K_X|$
containing $\fie\1(y)$, and $\fie\1(y)$ is contained in a fibre of
$\fie_{\om_{D'}}\colon D'\to\proj^1$, which contradicts Lemma~\ref{lem:hh}.
Hence $\fie\colon X\to Y$ is of degree 2 (possibly inseparable if $\chara
k=2$). In particular $2\mid9K^2$, so that $K^2$ is even and $K^2\ge4$; thus
$P_2\ge5$, and $\dim|2K_X-Z|\ge2$ for any cluster $Z$ of degree 2. By
changing $Z$ if necessary, we can assume that $\fie(Z)=y\in Y$ is a general
point, and is thus nonsingular. We have just shown that every fibre
$\fie\1(y)$ has degree exactly 2, so that $\fie$ is flat by
Lemma~\ref{lem:flat}; it is easy to see that this implies that $Y$ is normal.

Now for any $D\in|2K_X-Z|$, the image $\fie(D)=\Ga_D\subset Y$ is a curve
through $y=\fie(Z)$ isomorphic to $\proj^1$, and $\deg\fie_{|D}=\deg\fie=2$.
By Lemma~\ref{lem:push} the $\Ga_D\subset Y$ are linearly equivalent, so that
they are all contained in a linear system. This contradicts
Lemma~\ref{lem:sing}: in any linear system of curves through $y$, curves
singular at $y$ form a linear subsystem of codimension $\le2$, whereas the
$\Ga_D$ for $D\in|2K_X-Z|$ form an algebraic subfamily of nonsingular curves
depending with a complete parameter space of dimension $\ge2$ made up of
curves isomorphic to $\proj^1$. \QED \end{pfof}

 \begin{REM} Here we have assumed that $\fie(Z)=y\in Y$ is a general
point only for simplicity (see Lemma~\ref{lem:sing}).
 \end{REM}

 \begin{pfof}{Theorem~\ref{th:tri}, Case~(b)} Let $Z$ be a cluster of degree
2 on $X$ and $x\in Z$ a reduced point; that is, $Z$ is either a pair
$(x,y)$ of distinct points, or a point $x$ plus a tangent vector $y$ at $x$.
We assume that $|3K_X|$ does not separate $Z$, and gather together a number
of deductions concerning the curves
 \begin{equation}
 C_L\in|K_X+L|\quad\text{and}\quad
D_L\in|2K_X-L|\quad\text{for all}\quad L\in\Pico X,
 \nonumber
 \end{equation}
arriving eventually at a contradiction.

 \subparagraph{Step A} $h^0(K_X+L)\ge1$ for all $L\in\Pico X$. In fact if
$L\ne0$ then $h^2(K_X+L)=0$, and hence $h^0(K_X+L)\ge\chi(K_X)\ge1$.

 \subparagraph{Step B} $Z\not\subset C_L$ for all $L\in\Pico X$ and all
$C_L\in|K_X+L|$. Indeed
 \begin{equation}
H^0(X,\Oh_X(3K))\to H^0(C_L,\Oh_{C_L}(3K_X))
 \nonumber
 \end{equation}
is onto by the assumption $H^1(\Oh_X(2K_X-L))=0$, and $\Oh_{C_L}(3K_X)$ very
ample follows from Theorem~\ref{th:curve} exactly as in \S\ref{sec:pluri},
Step~II. Therefore if $Z\subset C_L$ then $|3K_X|$ separates $Z$, which we
are assuming is not the case.

 \subparagraph{Step C} For general $L\in\Pico X$ and all $C_L\in|K_X+L|$ we
have $x\in C_L$.

First of all, since $\dim\Pico X\ge1$, there is an $L\in\Pico X$ and a curve
$C_L\in|K_X+L|$ containing $x$, and $C_L$ does not contain $Z$ by Step~B.
Now if $L_1,L_2\in\Pico X$ is a general solution of $L+L_1+L_2=0$, and $x\notin
C_{L_1}$, $x\notin C_{L_2}$, then $C_L+C_{L_1}+C_{L_2}$ separates $x$ and
$Z$, a contradiction.

 \subparagraph{Step D} $h^0(K_X+L)=1$ and $H^1(K_X+L)=0$ for general $L\in\Pico
X$.

By Step~C, every $s\in H^0(K_X+L)$ vanishes at $x$. If $h^0(K_X+L)\ge2$ then
some nonzero section would vanish also at $y$. The statement about $H^1$
follows from RR:
 \begin{equation}
1=h^0(\Oh_X(K_X+L))\ge\chi(\Oh_X(K_X+L))=\chi(\Oh_X)\ge1.
 \nonumber
 \end{equation}

 \subparagraph{Step E} $x\in\Bs|2K_X-L|$ for general $L\in\Pico X$. For if
$D_L\in|2K_X-L|$ does not contain $x$ then $D_L+C_L$ separates $x$ from
$Z$ (since by Step~B already $C_L$ separates them).

 \subparagraph{Step F} For general $L,L_1\in\Pico X$, the point $x$ is a base
point of the linear system $\bigl|(2K_X-L_1)_{\textstyle{|C_L}}\bigr|$ on
$C_L$, and hence
 \begin{equation}
H^1(m_x\Oh_{C_L}(2K_X-L_1))\ne0.
 \nonumber
 \end{equation}

This follows from $x\in\Bs|2K_X-L_1|$ because by Step~D, restriction from $X$
maps onto $H^0(\Oh_{C_L}(2K_X-L_1))$.

 \subparagraph{Step G} We now observe that Step~B implies that $x$ is a
singular point of $C_L$. If $x\in\Sing X$ then it is automatically singular on
$C_L$. On the other hand, if $x$ is nonsingular on $X$ and on $C_L$, consider
the blowup $\si\colon X_1\to X$ of $x$ and the algebraic system
$C'_L=\si^*C_L-E$, where $E$ is the exceptional divisor. Let $y\in X_1$ be the
point corresponding either to the other point or to the tangent vector of the
cluster $Z$. Since the curves $C'_L$ move in a positive dimensional system,
there is a curve $C'_L$ through $y$, and therefore a curve $C_L$ containing
$Z$, contradicting Step~B.

 \subparagraph{Step H} For general elements $L,L_2\in\Pico X$, there is an
isomorphism $m_x\Oh_{C_L}(L_2)\iso m_x$.

This follows as usual by automatic adjunction (Lemma~\ref{lem:adj}) applied to
the conclusion $H^1(m_x\Oh_{C_L}(2K_X-L_1))\ne0$ of Step~F, where
$L_1=-L-L_2$. We first get a nonzero homomorphism
 \begin{equation}
m_x\Oh_{C_L}(2K_X-L_1)\to\om_{C_L}=\Oh_{C_L}(2K_X+L),
 \nonumber
 \end{equation}
that is, a map $m_x\Oh_{C_L}(L_2)\to\Oh_{C_L}$; since $C_L$ is 2-connected
this must be an inclusion, and the image is the ideal of a point $m_z$. But
$x$ is a singular point of $C_L$ (by Step~G), and thus $x=z$.

 \subparagraph{Step I} Let $\si\colon C'\to C=C_L$ be the blowup at $x$.
Step~H implies that $L_2'=\si^*L_2$ is trivial on $C'$ for every general $L_2$,
and hence for every $L_2\in\Pico X$ (by the group law). We derive a
contradiction from this. Consider the diagram
 \begin{equation}
 \renewcommand\arraystretch{1.5}
 \begin{array}{rcl}
\Pico X @>\res_C>> & \Pico C & @>\si^*>> \Pico(C')\\
& \uparrow \\
& G
 \end{array}
 \nonumber
 \end{equation}
where $G$ is the kernel of $\si^*$. Now the key point (exactly as in Ramanujam
and Francia vanishing) is that $G$ is an affine group scheme. Since the
composite $\si^*\circ\res_C$ is zero, $\Pico X$ maps to $G$. Since $\Pico
X$ is complete $\res_C$ is the constant morphism to zero. But this is
obviously nonsense: for example, since $H^1(\Oh_X(2K_X+L))=0$ for all
$L\in\Pico X$, the exact sequence
 \begin{equation}
0\to\Oh_X(-K_X-L+N)\to\Oh_X(N)\to\Oh_C(N)\to0
 \nonumber
 \end{equation}
is exact on global sections if $L\ne N$. Thus $H^0(\Oh_C(N))=0$ and
the restriction of $N$ to $C$ is nontrivial. \QED \end{pfof}

 \section{The bicanonical map}\label{sec:bi}
 \subsection*{Preliminaries and the proof of Theorem~\ref{th:bi}, (a) and (c)}
This section proves Theorem~\ref{th:bi}. We start by remarking that $|2K_X|$ is
free. Indeed, for any $C\in|K_X|$, the restriction $\Oh_X(2K_X)\to\Oh_C(K_C)$
is surjective on $H^0$, and $|K_C|$ is free by Theorem~\ref{th:free}. For a
cluster $Z$ of degree 2 in $X$, note the following obvious facts:
 \begin{enumerate}
 \renewcommand{\labelenumi}{(\roman{enumi})}
 \item If $Z$ is contracted by $|2K_X|$ then $|K_X|$ does not separate
$Z$; thus
 \begin{equation}
h^0(\sI_Z\Oh_X(K_X))\ge p_g-1\quad\text{or}\quad \dim|K_X-Z|\ge p_g-2.
 \nonumber
 \end{equation}
 \item If $|2K_X|$ contracts $Z$ then so does $|K_C|$ for any curve
$C\in|K_X-Z|$.
 \end{enumerate}

 \begin{pfof}{Theorem~\ref{th:bi}, (a)} We suppose that every curve
$C\in|K_X-Z|$ is \hbox{3-connected}, and derive a contradiction from the
assumption that $|2K_X|$ contracts $Z$. By Theorem~\ref{th:hh}, every
$C\in|K_X-Z|$ is honestly hyperelliptic. As in the proof of
Theorem~\ref{th:tri}, Case~(a), it follows that $\fie_{2K}\colon X\to Y$ has
degree 2, and maps every $C\in|K_X-Z|$ as a double cover of a curve
$\Ga_C\subset Y$ isomorphic to $\proj^1$. Then $\Ga_C$ for $C\in|K_X-Z|$
form an algebraic subfamily of a linear system of curves through
$y=\fie_{2K}(Z)$, with a complete parameter space of dimension $\ge2$. As
before, this contradicts Lemma~\ref{lem:sing} (but $y\in Y$ may now be
singular). \QED
 \end{pfof}

 \begin{DEF}\label{def:Frc} Let $X$ be a projective surface with at worst Du
Val singularities and with $K_X$ nef. A {\em Francia curve} or {\em Francia
cycle} is an effective Weil divisor $B$ on $X$ satisfying 
 \begin{equation}
 K_XB=p_aB=1\text{ or }2.
 \nonumber
 \end{equation}
If $K_X$ is ample and $B$ is Gorenstein (for example if $B$ is a Cartier
divisor), it is clearly either an irreducible curve of genus 1, or a
numerically 2-connected curve of arithmetic genus $p_a=2$. It would be
interesting to know if $B$ is necessarily Gorenstein.
 \end{DEF}

 \begin{pfof}{Theorem~\ref{th:bi}, (b) $\implies$ (c)} The argument is standard
and we omit some details. Suppose that the 2-canonical map
$\fie=\fie_{2K}\colon X\to Y$ is not birational. Every point $x\in X$ is
contained in a cluster $Z$ of degree 2 contracted by $\fie$; we choose
$x\in\NonSing X$. Theorem~\ref{th:bi}, (b) gives a Francia curve $B_0\subset X$
through $Z$. Write $S\to X$ for the minimal nonsingular model of $X$ and
$B=\Bh_0$ for the hat transform of $B_0$ (as in the proof of
Lemma~\ref{lem:n-conn}). Then by Claim~\ref{cla:bhat}, $B$ is also a Francia
cycle on $S$, that is, $1\le K_SB=p_aB\le2$. An easy argument in quadratic
forms shows that there are at most finitely many effective divisors $B\subset
S$ with $K_SB=1$ and $B^2=-1$ (compare \cite{Bo}, pp.~191--192 or \cite{BPV},
p.~224). Therefore every general point of $S$ is contained in a curve $B$ with
$K_SB=p_aB=2$, and hence $B^2=0$. Now the same argument in quadratic forms
shows that divisors with $K_SB=2$ and $B^2=0$ belong to finitely many
numerical equivalence classes, so one class must contain an algebraic family
of curves. This gives a genus 2 pencil on $S$, and therefore also on $X$. \QED
\end{pfof}

We use the following obvious lemma at several points in what follows.

 \begin{LEM}[Dimension lemma]\label{lem:dim} Let $\eta\subset X$ be a cluster
of degree $d$ which is contracted by $|2K_X|$, and $C\in|K_X|$ a curve
containing $\eta$. Then
 \begin{equation}
h^1(\sI_\eta\Oh_C(K_C))=\dim\Hom(\sI_\eta,\Oh_C)=d.
 \nonumber
 \end{equation}
In particular, for any $x\in C$, we have
 \begin{equation}
h^1(m_x^2\Oh_C(K_C))=\dim\Hom(m_x^2\Oh_C,\Oh_C)=1+\dim T_{\fie,x}\le4,
 \nonumber
 \end{equation}
where $T_{\fie,x}$ is the Zariski tangent space to the scheme theoretic fibre
of $\fie_{2K_X}$ through $x$. \end{LEM}

 \begin{pf} Since $|K_C|$ is free and contracts $\eta$, the evaluation map
$H^0(\Oh_C(K_C))\to\Oh_\eta(K_C)=k^d$ has rank 1, so that
$h^1(\sI_\eta\Oh_C(K_C))=d$ comes from the exact sequence
 \begin{equation}
 \renewcommand\arraystretch{1.3}
 \begin{array}{l}
0\to H^0(\sI_\eta\Oh_C(K_C))\to H^0(\Oh_C(K_C))\to k^d \\
\hphantom{0}\to H^1(\sI_\eta\Oh_C(K_C))\to H^1(\Oh_C(K_C))=k.
 \end{array}
 \nonumber
 \end{equation}
As usual, Serre duality gives
 \begin{equation}
\Hom(\sI_\eta,\Oh_C)=\Hom(\sI_\eta\Oh_C(K_C),\om_C)\dual H^1(\sI_\eta\Oh_C(K_C)).
 \nonumber
 \end{equation}
We obtain the last part by taking $\eta$ to be the intersection of the scheme
theoretic fibre $\fie\1(\fie(x))$ with the subscheme $V(m_x^2)\subset C$
corresponding to the tangent space. \QED \end{pf}

 \subsection*{Case division and plan of proof of~(b)}

Throughout this section, $Z$ is a cluster of degree 2, and we argue by
restricting to a curve $C\in|K_X-Z|$, usually imposing singularities on $C$ at
a point $x\in Z$. As usual, the assumption that $Z$ is contracted by $K_C$
gives a homomorphism $\sI_Z\to\Oh_C$ linearly independent of the identity
inclusion. By passing to a suitable linear combination $s'=s+\la\id$ if
necessary, we assume that $s\in\Hom(\sI_Z,\Oh_C)$ is injective, and hence
$s(\sI_Z)=\sI_{Z'}$ for some cluster $Z'$ of degree 2; the family of clusters
$Z'$ as $s$ runs through injective elements $s\in\Hom(\sI_Z,\Oh_C)$ is an
analog of a $g^1_2$ on $C$.

The argument is modelled on the proof of Theorem~\ref{th:hh}. As there, we use
different arguments depending on how $Z$ and $Z'$ intersect, or, to put it
another way, how $Z'$ moves as $s$ runs through injective elements
$s\in\Hom(\sI_Z,\Oh_C)$. (In other words, how the $g^1_2$ corresponding to
$\Hom(\sI_Z,\Oh_C)$ breaks up into a ``base locus'' plus a ``moving part''.)
Let $s\in\Hom(\sI_Z,\Oh_C)$ be a general element, and $\sI_{Z'}=s(\sI_Z)$.
Logically, there are 4 cases for $Z$ and $Z'$.

 \begin{enumerate}
 \item $\Supp Z\cap\Supp Z'=\emptyset$.
 \item $\Supp Z\cap\Supp Z'\ne\emptyset$, but $\Supp Z\ne\Supp Z'$.
 \item $Z=Z'$.
 \item $Z\ne Z'$ are nonreduced clusters supported at the same point $x\in
X$.
 \end{enumerate}

In Case~2, $|Z|$ has a fixed point plus a moving point; as we see in
Lemma~\ref{lem:case2}, this contradicts $K_X$ ample. In Case~1, $|Z|$ is a
free $g^1_2$, and the isomorphism $\sI_Z\iso\sI_{Z'}$ with $\Supp Z\cap\Supp
Z'=\emptyset$ implies that $\sI_Z$ is locally free, so that $Z$ is a Cartier
divisor on $C$. If $p_g\ge4$, it turns out that we can choose $C$ to be
``sufficiently singular'' at a point $x\in Z$ so that $Z\subset C$ is not
Cartier, and Case~1 is excluded for such $C$ (see Lemma~\ref{lem:notCt}).

In Cases~3--4, when the support of $Z$ does not move, we must find a map
$s'\colon\sI_Z\to\Oh_C$ vanishing on a ``fairly large'' portion of $C$, so that
its scheme theoretic support $B\subset C$ is ``fairly small''. The key idea is
to look for $s'$ as a nilpotent or idempotent (see Lemma~\ref{lem:pot} and
Corollary~\ref{cor:Art}). The assumption of Case~3 is
$\Hom(\sI_Z,\Oh_C)=\End(\sI_Z)$, which is a 2-dimensional Artinian algebra;
this makes it is rather easy to find a nilpotent or idempotent element, and to
prove Theorem~\ref{th:bi}, (b).

In Case~4, $Z'$ is $x$ plus a tangent vector $y$ which moves in $T_{C,x}$ as
$s\in\Hom(\sI_Z,\Oh_C)$ runs through injective elements; this is an {\em
infinitesimal} $g^1_2$, an interesting geometric phenomenon in its own right
(see Remark~\ref{rem:g23} and the proof of Proposition~\ref{pro:m2}, Step~6 for
more details). The key point in this case is to prove that the extra
homomorphism $s\colon\sI_Z\to\Oh_C$ takes $m_x^2$ to itself, so that
$\End(m_x^2)$ is a nontrivial Artinian algebra; see Proposition~\ref{pro:m2}.

 \begin{REM}\label{rem:g23} In Case~4, reversing the usual argument proves
that $\fie_{K_C}$ also contracts $Z'$, and so it contracts a cluster $\eta$ of
degree $\ge3$ contained in the first order tangent scheme $V(m_x^2)\subset C$.
If $C$ is numerically 3-connected, this is of course impossible by
Theorem~\ref{th:hh}. In this case, $\Hom(\sI_\eta,\Oh_C)$ is a certain analog
of a $g^2_3$ or $g^3_4$ on $C$.

Case~4 certainly happens on abstract numerically 2-connected Gorenstein
curves, and more generally, the analog of a $g^{m-1}_m$. Example: let $C_i$ for
$i=1,\dots,m$ be nonhyperelliptic curves of genus $g_i\ge3$ with marked points
$x_i\in C_i$, and assemble the $C_i$ into a curve $C=\bigcup C_i$ by glueing
together all the $x_i$ to one point $x$, at which the tangent directions are
subject to a single nondegenerate linear relation, so that the singularity
$x\in C$ is analytically equivalent to the cone over a frame of reference
$\{P_1,P_2,\dots,P_m\}$ in $\proj^{m-2}$. Then $C$ is Gorenstein and
$K_C$ restricted to each $C_i$ is $K_{C_i}+2x_i$ (see \cite{Ca1},
Proposition~1.18, (b), p.~64, or \cite{Re}, Theorem~3.7), so that $|K_C|$
contracts the whole $(m-1)$-dimensional tangent space $T_{C,x}$ to a point.

A cluster $Z$ of degree 2 supported at $x$ corresponds to a point
$Q\in\proj^{m-2}=\proj(T_{C,x})$. Since $Z$ is contracted by $K_C$ (together
with the whole tangent space), by our usual argument, the group
$\Hom(\sI_Z,\Oh_C)$ is 2-dimensional and a general $s\colon\sI_Z\to\Oh_C$
has image $\sI_{Z'}$ where $Z'$ is a moving cluster of degree 2 at $x$,
corresponding to a moving point $Q'\in\proj^{m-2}$. It is an amusing exercise
to see that if $Q$ is linearly in general position with respect to the frame
of reference $\{P_1,P_2,\dots,P_m\}$ then $Q'$ moves around the unique
rational normal curve of degree $m-2$ passing through
$\{P_1,P_2,\dots,P_m,Q\}$. On the other hand, if $Z$ is in the tangent cone
to $C$ (say, tangent to the branch $C_1$), then $\sI_Z$ is not isomorphic to
any other cluster of degree 2, so that $\Hom(\sI_Z,\Oh_C)=\End(\sI_Z)$;
this has 2 idempotents vanishing on $C_1$ and on $C_2+\cdots+C_m$.

The following easy exercises may help to clarify things for the reader:
 \begin{enumerate}
 \item Let $x\in C$ be an ordinary triple point of a plane curve, say defined
by an equation $f(u,v)=u^3+v^3+$ higher order terms; then for general $\la$,
the ideals $(u+\la v,v^2)$ in $\Oh_{C,x}$ are all locally isomorphic. [Hint:
Multiply by the rational function $(u+\mu v)/(u+\la v)$.]
 \item If $C$ is the planar curve defined by $vw=v^3+w^3$ then $m_x=(v,w)$ is
locally isomorphic to $\sI_Z=(v,w^2)$ and to $\sI_{Z'}=(v^2,w)$.
 \item If $C$ is the planar curve locally defined by $v^2=w^3$ then $m_x=(v,w)$
is locally isomorphic to $\sI_Z=(v,w^2)$.
 \end{enumerate}
(Compare the proof of Proposition~\ref{pro:m2}, Step~6.)
 \end{REM}

 \begin{LEM}\label{lem:case2} Case~2 is impossible. \end{LEM}

 \begin{pf} Since $x\in Z\cap Z'$ and $\Supp Z\ne\Supp Z'$, we can interchange
$Z$ and $Z'$ if necessary and assume that $Z'=\{x,y\}$ with $x\ne y$. Consider
the inclusion $s\colon\sI_Z\into\Oh_C$ with image $s(\sI_Z)=\sI_{Z'}=m_xm_y$
and the identity inclusion. One of these vanishes at $y$ and the other
doesn't, so their restrictions to a component $\Ga$ containing $y$ are linearly
independent on $\Ga$, and, as in Claim~\ref{cla:mov_y}, for any general point
$y'\in\Ga$, some linear combination $s'=s+\la\id$ defines an isomorphism
$s'\colon\sI_Z\iso m_xm_{y'}$. Reversing our usual argument shows that $x$ and $y'$ are
contracted to the same point by $|K_C|$ or $|2K_X|$, so that the free linear
system $|2K_X|$ contracts $\Ga$ to a point. This contradicts $K_X$ ample. \QED
\end{pf}

 \subsection*{Clusters on singular curves}
Our immediate aim is to exclude Case~1, but at the same time we introduce some
ideas and notation used throughout the rest of this section. Choose a point
$x\in Z$. Since $X$ has at worst hypersurface singularities and $C$ is a
Cartier divisor in $X$, it is a local complete intersection, that is, locally
defined by $F=G=0$. (Of course, $X$ may be nonsingular.) We think of $x\in
Z\subset C\subset X\subset\aff^3$ as local, and write $\Oh_{\aff^3}$, $\Oh_C$,
etc.\ for the local rings at $x$. We take local coordinates $u,v,w$ in $\aff^3$
so that $Z$ is defined by $u=v=w=0$ in the reduced case, or $u=v=w^2=0$
otherwise.

 \begin{LEM}\label{lem:notCt}
 \begin{enumerate}
 \renewcommand{\labelenumi}{{\rm(\arabic{enumi})}}
 \item The quotient $\sI_{\aff^3,Z}/m_{\aff^3,x}\sI_{\aff^3,Z}$ is a
$3$-dimensional vector space, and $Z\subset C$ is a Cartier divisor at $x$ if
and only if $F,G$ map to linearly independent elements of it.
 \item Suppose that $p_g\ge4$ and $Z$ is contracted by $|2K_X|$. Then the curve
$C\in|K_X-Z|$ can be chosen such that $Z$ is not a Cartier divisor. For this
$C$,  Case~1 is excluded. 
 \end{enumerate}
 \end{LEM}

 \begin{pf} (1) says that a minimal set of generators of the ideal
$\sI_{\aff^3,Z}$ consists of 3 elements, which is obvious because
$\sI_{\aff^3,Z}$ is locally generated at $x\in Z$ by the regular sequence
$(u,v,w)$ or $(u,v,w^2)$. Now $Z$ is a Cartier divisor on $C$ if and only if
$\sI_{C,Z}$ is generated by 1 element, that is, $F$ and $G$ provide two of the
minimal generators of $\sI_{\aff^3,Z}$. This proves (1).

For (2), suppose that $F=0$ is the local equation of $X\subset\aff^3$. If $F\in
m_{\aff^3,x}\sI_{\aff^3,Z}$ then by (1), $Z$ is not a Cartier divisor on any
curve $C\in|K_X-Z|$. Suppose then that $F\notin m_{\aff^3,x}\sI_{\aff^3,Z}$, so
that $F$ provides one of the minimal generators of $\sI_{\aff^3,Z}$. Then the
ideal $\sI_{X,Z}$ of $Z\subset X$ is generated by 2 elements, in other words,
$\dim_k\sI_{X,Z}/m_{X,x}\sI_{X,Z}=2$. Therefore
 \begin{equation}
h^0(m_x\sI_Z\Oh_X(K_X))\ge h^0(\sI_Z\Oh_X(K_X))-2\ge p_g-3\ge1
 \nonumber
 \end{equation}
(by remark (i) at the beginning of this section). Thus we can find a curve
$C\in|K_X-Z|$ whose local equation at $x$ is $g\in m_{X,x}\sI_{X,Z}$. Then
$g$ has a local lift $G\in m_{\aff^3,x}\sI_{\aff^3,Z}$, so that (1) applies to
$C$. \QED \end{pf}

 \begin{REM}\label{rem:geom} The same argument can be expressed more
geometrically. If $Z$ contains $x$ as a reduced point, that is,
$\sI_{\aff^3,Z}=m_x$, then $x\in C$ is Cartier if and only if $C$ defined by
$(F,G)$ is nonsingular at $x$, that is, $F,G$ map to linearly independent
elements of $m_x/m_x^2$.

To interpret the nonreduced case $\sI_{\aff^3,Z}=(u,v,w^2)$, note that
 \begin{equation}
F\notin m_{\aff^3,x}\sI_{\aff^3,Z} \iff F=Pu+Qv+Rw^2
\quad\text{with one of $P,Q,R\notin m_x$.}
 \nonumber
 \end{equation}
In other words, the surface $Y$ locally defined by $F=0$ is either nonsingular
at $x$, or has a double point with $Z$ not in the tangent cone. In the opposite
case $F\in m_{\aff^3,x}\sI_{\aff^3,Z}$, it is easy to see that $x\in C$ is
either a complete intersection defined by two singular hypersurfaces, so has
3-dimensional tangent space $T_{C,x}$, or is a planar curve, which is either a
double point with $Z$ in the tangent cone, or a point of multiplicity $\ge3$.
 \end{REM}

 \subsection*{The nilpotent--idempotent lemma}
Our proof of Theorem~\ref{th:bi}, (b) in Cases~3--4 is based on the following
result. Note first that $\Hom(\sI_Z,\Oh_C)\subset H^0(C\setminus\Supp
Z,\Oh_C)$, and the latter is a ring. (We usually write $\sI_Z$ for
$\sI_{C,Z}$ in what follows.) In other words, maps $\sI_Z\to\Oh_C$ can be
viewed as rational sections of $\Oh_C$ that are regular outside $\Supp Z$, so
that it is meaningful to multiply them (the product is again a rational
section of $\Oh_C$ that is regular away from $Z$).

 \begin{LEM}\label{lem:pot} Assume that $K_X^2\ge10$, and let $C\in|K_X-Z|$.
Suppose that $s\colon\sI_Z\to\Oh_C$ is a nonzero homomorphism which is either
nilpotent with $s^4=0$, or a nontrivial idempotent with $s(1-s)=0$. Then the
scheme theoretic support of $s$ (respectively, in the idempotent case, either
$s$ or $1-s$) is a Francia curve $B$, and $\sI_Z\Oh_B(2K_X)\iso\om_B$.

More generally, suppose that $s_i\colon\sI_Z\to\Oh_C$ for $i=1,\dots,4$ are
nonzero homomorphisms such that $s_1s_2s_3s_4=0$. Then one of the $s_i$ has
scheme theoretic support a Francia curve $B_i$ with
$\sI_Z\Oh_{B_i}(2K_X)\iso\om_{B_i}$.
 \end{LEM}

The final part is more general, because we allow some $s_i=\id$, or some of
the $s_i$ to coincide. Notice that $\Oh_C$ has no sections supported at
finitely many points, so we need only check the conditions $s^4=0$ etc.\
in each generic stalk of $\Oh_C$, that is, as rational functions on $C$.

 \begin{pf} If $s\colon\sI_Z\Oh_C(K_C)\to\om_C$ is not generically injective,
the factorisation provided by automatic adjunction (Lemma~\ref{lem:adj}) gives
a subcurve $B\subset C$ satisfying $\sI_Z\Oh_B(K_C)\iso\om_B$; we are in the
limiting case of numerically 2-connected. Write $C=A+B$ for the decomposition
of Weil divisors, so that $A$ is the divisor of zeros of $s$. Passing to the
minimal nonsingular model $S$ and taking the hat transform $\Bh$ as in
Lemma~\ref{lem:n-conn} and Claim~\ref{cla:bhat} gives a decomposition
$K_S\lineq f^*C=A_1+\Bh$ such that $A_1\Bh=2$.

Therefore by the Hodge algebraic index theorem, $A_1^2\Bh^2\le(A_1\Bh)^2=4$. If
both $A_1^2$, $\Bh^2\ge1$, it follows that $K_S^2\le9$, a contradiction, so
that either $A_1^2\le0$ or $\Bh^2\le0$. Then (because $K_S=A_1+\Bh$ and
$A_1\Bh=2$), either $K_XA=K_SA_1\le2$ or $K_XB=K_S\Bh\le2$. Suppose for the
moment that $K_S\Bh\le2$. Since $2p_a\Bh-2=\Bh^2+K_S\Bh$, it follows at once
that we are in one of the two cases
 \begin{equation}
\Bh^2=-1,K_S\Bh=1,p_a\Bh=1\quad\text{or}\quad
\Bh^2=0,K_S\Bh=2,p_a\Bh=2.
 \nonumber
 \end{equation}
But by Lemma~\ref{lem:n-conn} and Claim~\ref{cla:bhat} we have $K_S\Bh=K_XB$
and $p_a\Bh=p_aB$, so that $B$ is the required Francia curve.

It remains to get rid of the possibility that $K_XA=K_SA_1\le2$ in the
different cases. If $s$ is a nontrivial idempotent, we can swap $A\bij B$ by
$s\bij1-s$ if necessary, so that $K_XB\le2$. In the nilpotent case, since
$A$ equals the Weil divisor of zeros of $s$ and $s^4=0$, it follows that
$C\le4A$. Then $K_XA\le2$ would imply $K_X^2\le8$, a contradiction.

The last part is exactly the same: each $s_i$ (for $i=1,2,3,4$) is either
injective or has scheme theoretic support a subcurve $B_i\subset C$ with
$\sI_Z\Oh_{B_i}(K_C)\iso\om_{B_i}$, and divisor of zeros $A_i=C-B_i$. Since
$\prod s_i=0$ it follows that $C\le\sum A_i$. Now arguing as above gives that
one of $K_XA_i$ or $K_XB_i\le2$; if the first alternative holds for all $i$
then $K_X^2=K_XC\le \sum K_XA_i\le8$, a contradiction. This proves the lemma.
\QED \end{pf}

We apply Lemma~\ref{lem:pot} via a simple algebraic trick.

 \begin{COR}\label{cor:Art} If $A=\End_{\Oh_C}(\sI_{C,Z})$ is an Artinian
algebra of length $\ge2$ then it has a nontrivial idempotent or a nonzero
nilpotent with $s^2=0$. More generally, if\/ $\Hom(\sI_{C,Z},\Oh_C)$ is a
$2$-dimensional vector space contained in an Artinian algebra $A\subset
H^0(C\setminus\Supp Z,\Oh_C)$ of dimension $\le4$ then there exist nonzero
elements $s_1,\dots,s_4\in\Hom(\sI_{C,Z},\Oh_C)$ with zero product. Under
either assumption, Lemma~\ref{lem:pot} gives a Francia curve $B\subset C$
containing $Z$. \end{COR}

This completes the proof of Theorem~\ref{th:bi}, (b) in Case~3, since the case
assumption is that $s\colon\sI_Z\to\sI_Z\subset\Oh_C$, so that
$\Hom(\sI_Z,\Oh_C)=\End(\sI_Z)$ is a 2-dimensional Artinian algebra.

 \begin{pf} In the main case $\dim A=2$, this is completely trivial: if
$k\subset A$ is the constant subfield, any $s\in A\setminus k$ satisfies a
quadratic equation over $k$ of the form
 \begin{equation}
0=s^2+as+b=(s-\al)(s-\be).
 \nonumber
 \end{equation}
If $\al=\be$ then $s'=s-\al$ is nilpotent with $s'{}^2=0$; otherwise,
$s'=(s-\al)/(\al-\be)$ and $1-s'=(s-\be)/(\be-\al)$ are orthogonal idempotents.

More generally, an Artinian algebra is a product $A=A_1\times\cdots\times A_l$
with local Artinian rings $(A_i,n_i)$ as factors; the maximal ideals of $A$
are codimension 1 vector subspaces $m_i\subset A_i$ given by $n_1\times
A_2\times\cdots\times A_l$ (say). The projection to the factors (if $l\ge2$)
give nontrivial idempotents; if $l=1$ then $A$ is local, with nilpotent maximal
ideal. This proves the first part.

We now prove the more general statement: a 2-dimensional vector subspace
$V\subset A$ in an Artinian algebra has nonzero intersection with every
maximal ideal, say $s_i\in V\cap m_i$. If the local factors $(A_i,n_i)$ have
dimension $d_i$ then $n_i^{d_i}=0$, and the product $\prod s_i^{d_i}$ maps to
zero in each factor, so is zero in $A$. Taking $\sum d_i=\dim A\le4$ gives the
final part of the claim. \QED \end{pf}

 \subsection*{Proof in Case~4}
In the following proposition, $x\in C\subset\aff^3$ is a {\em local} curve
which is a local complete intersection at $x$. We choose local coordinates
$u,v,w$ on $\aff^3$ so that $\sI_{\aff^3,Z}\subset\Oh_{\aff^3}$ is generated at
$x$ by the regular sequence $u,v,w^2$. As before, we write $\Oh_C$ for the
local ring $\Oh_{C,x}$ and $\sI_Z=\sI_{C,Z}$ for the $\Oh_C$ module obtained
as the stalk at $x$ of the corresponding ideal sheaf. (Thus the statement of
the proposition only concerns homomorphisms $s\colon\sI_Z\to\Oh_C$ of
modules over the local ring $\Oh_C$.)

 \begin{PROP}\label{pro:m2} Let $Z\subset C$ be a cluster of degree $2$
supported at $x$. We assume
 \begin{enumerate}
 \renewcommand{\labelenumi}{\rm(\roman{enumi})}
 \item $Z$ is not a Cartier divisor on $C$;
 \item there exists a homomorphism $s_0\colon\sI_Z\to\Oh_C$ such that for
general $\la\in k$, $s_0+\la\id$ defines an isomorphism $\sI_Z\iso\sI_{Z_\la}$
with $Z_\la$ a cluster of degree $2$ supported at $x$, and $Z_0\ne Z$.
 \end{enumerate}
Then any homomorphism $s\colon\sI_Z\to\Oh_C$ takes $m_{C,x}^2$ to
$m_{C,x}^2$, that is,
 \begin{equation}
 \renewcommand\arraystretch{1.5}
 \begin{matrix}
\sI_Z & @>{\quad s\quad}>> & \Oh_C\\
\bigcup && \bigcup\\
m_x^2 & @>{\hphantom{\quad s\quad}}>> & m_x^2
 \end{matrix}
 \nonumber
 \end{equation}
 \end{PROP}

 \begin{pfof}{Theorem~\ref{th:bi}, (b) in Case~4} We apply the proposition to
the {\em global} homomorphism $s\colon\sI_Z\to\Oh_C$, using the assumption
of Case~4. We get
 \begin{equation}
\Hom(\sI_Z,\Oh_C) \subset \End(m_x^2) \subset \Hom(m_x^2,\Oh_C).
 \nonumber
 \end{equation}
Now Lemma~\ref{lem:dim} gives $\dim\Hom(\sI_Z,\Oh_C)=2$ and
$\dim\Hom(m_x^2,\Oh_C)\le4$; but $A=\End(m_x^2)$ is a subring of
$H^0(C\setminus\Supp Z,\Oh_C)$, so that Corollary~\ref{cor:Art} gives the
result. \QED \end{pfof}

 \begin{pfof}{Proposition~\ref{pro:m2}, Step 1} If $s\in\Hom(\sI_Z,\Oh_C)$
is any element then $s(\sI_Z)\subset m_x$; for otherwise $s$ would be an
isomorphism
$\sI_Z\iso\Oh_C$ near $x$, contradicting the assumption that $Z\subset C$ is
not Cartier.

 \subparagraph{Step 2} Note that $m_x^2\subset\sI_Z$, so that we can
restrict $s\colon\sI_Z\to\Oh_C$ to $m_x^2$. Also, $m_x\sI_Z\subset m_x^2$,
and obviously $s(\sI_Z)\subset m_x$ implies that $s(m_x\sI_Z)\subset m_x^2$.

 \subparagraph{Step 3} It is enough to prove that $s(w^2)\in m_x^2$. Indeed,
 \begin{equation}
 m_x\sI_Z=(u,v,w)\cdot(u,v,w^2)=(u^2,uv,v^2,uw,vw,w^3),
 \nonumber
 \end{equation}
so that
 \begin{equation}
 m_x^2=(u,v,w)^2=(u^2,uv,v^2,uw,vw,w^2)=m_x\sI_Z+\Oh_Cw^2\subset\Oh_C.
 \nonumber
 \end{equation}

 \subparagraph{Step 4} Since $C$ is a local complete intersection,
$\sI_{\aff^3,C}=(F,G)$, where $F,G\in\Oh_{\aff^3}$ is a regular sequence. Now
$Z\subset C$ gives $F,G\in\sI_{\aff^3,Z}$, so that
 \begin{equation}
 \begin{aligned}
 F&=Pu+Qv+Rw^2,\\ G&=P'u+Q'v+R'w^2,
 \end{aligned}
 \quad\text{with}\quad P,Q,R,P',Q',R'\in\Oh_{\aff^3}.
 \end{equation}

The set of local homomorphisms $\sI_Z\to\Oh_C$ is a module over $\Oh_C$;
this is the stalk at $x$ of the sheaf $\sHom$. For the moment, we take on
trust the following general fact (see Appendix to \S\ref{sec:bi} for a
discussion and a detailed proof.)

 \begin{CLA}\label{cla:PQ-PQ} The $\Oh_C$ module
$\sHom_{\Oh_C}(\sI_Z,\Oh_C)$ is generated by two elements, the identity
inclusion $\id\colon\sI_{C,Z}\into\Oh_C$ and the map
$t\colon\sI_{C,Z}\to\Oh_C$ determined by the minors of the $2\times3$ matrix
of coefficients of $F,G$:
 \begin{equation}
t(u)=QR'-RQ',\quad t(v)=-PR'+RP',\quad t(w^2)=PQ'-QP'.
 \end{equation}
 \end{CLA}

 \subparagraph{Step 5} According to Steps~3--4, to prove
Proposition~\ref{pro:m2}, we need only prove that $PQ'-QP'\in
m_{\aff^3,x}^2$. We are home if all four of $P,Q,P',Q'\in m_x$. Thus in what
follows, we assume (say) that $P'\notin m_x$. Then $P'$ is a unit, and $G=0$
defines a nonsingular surface $Y$ containing $C$. Dividing by $P'$, we can
rewrite $G$ in the form $u=-(Q'/P')v-(R'/P')w^2$. Then subtracting a multiple
of this relation from $F$ gives $f=qv+rw^2$ as the local equation of $C\subset
Y$ (where $q=Q-PQ'/P'$ and $r=R-PR'/P'$).

Therefore it only remains to prove that if $C$ is the planar curve defined by
$f=qv+rw^2$, the two assumptions of Proposition~\ref{pro:m2} imply that
$q\in m_{Y,x}^2$. As in Lemma~\ref{lem:notCt}, assumption (i) implies that
$q,r\in m_{Y,x}$, so that $q\in m_{Y,x}^2$ is equivalent to saying that $x\in
C\subset Y$ has multiplicity $\ge3$

 \subparagraph{Step 6} Consider the linear terms of the given isomorphism
$s_0\colon\sI_Z\to\sI_{Z_0}$: 
 \begin{equation}
 s_0(v)=av+bw\mod m_{Y,x}^2,\quad s_0(w^2)=cv+dw\mod m_{Y,x}^2.
 \nonumber
 \end{equation}
Because $Z_0\ne Z$, it follows that $(b,d)\ne(0,0)$. However, if $b=0$ and
$d\ne0$, then for general $\la$, the two generators of
$\sI_{Z_\la}=(s_0(v)+\la v,s(w^2)+\la w^2)$ would have linearly independent
linear terms, so that $\sI_{Z_\la}=m_{C,x}$. This contradicts assumption (ii).
Therefore $b\ne0$, and $\sI_{Z_\la}$ has a generator with the {\em variable}
linear term $(a+\la)v+bw$. It follows that $Z_\la$ runs linearly around the
tangent space to $x$ in $C$.

Now we claim that $x\in C\subset Y$ is a planar curve singularity of
multiplicity $\ge3$. Indeed, the isomorphism $\sI_Z\iso\sI_{Z_\la}$ implies
that $Z_\la\subset C$ cannot be a Cartier divisor; but if $x\in C\subset Y$
were a double point, this would restrict $Z_\la$ to be in the tangent cone,
contradicting what we have just proved. This completes the proof of
Proposition~\ref{pro:m2}. \QED \end{pfof}

 \subsection*{Appendix: Proof of Claim~\ref{cla:PQ-PQ}}
We start by slightly generalising the set-up: let $\Oh_{\aff}$ be a local
ring, assumed to be regular (for simplicity only), and $x,y,z$ a regular
sequence generating a codimension 3 complete intersection ideal
$\sI_Z=(x,y,z)$. Consider a regular sequence $F,G\in\sI_Z$. Note that
 \begin{equation}
F=Px+Qy+Rz\quad\text{and}\quad G=P'x+Q'y+R'z
 \nonumber
 \end{equation}
for some $P,\dots,R'\in\Oh_{\aff}$. Write $\Oh_C=\Oh_{\aff}/(F,G)$ and
$\sI_{C,Z}=\sI_Z\Oh_C=(x,y,z)\subset\Oh_C$. (In the application,
$Z\subset\aff=\aff^3$ was a nonreduced cluster defined by $(x,y,z)=(u,v,w^2)$
and $C\subset\aff^3$ a complete intersection curve through $Z$.)

 \begin{LEM}
 \begin{enumerate}
 \renewcommand{\labelenumi}{\rm(\arabic{enumi})}
 \item A presentation of $\sI_{C,Z}$ over $\Oh_C$ is given by
 \begin{equation}
 \Oh_C^{\oplus5}@>M>>\Oh_C^{\oplus3}
 @>\left(\begin{matrix} x\\y\\z\end{matrix}\right)>>\sI_{C,Z}\to0,
 \quad\text{where}\quad
 \renewcommand\arraystretch{1.2}
 M=\left(\matrix
P&Q&R\\
P'&Q'&R'\\
0&z&-y\\
-z&0&x\\
y&-x&0
 \endmatrix\right).
 \nonumber
 \end{equation}

 \item $\sHom(\sI_{C,Z},\Oh_C)$ is generated over $\Oh_C$ by the two
elements $\id$ and $t$, where
 \begin{equation}
t\colon\left(\matrix x\\y\\z\endmatrix\right)\mapsto
\left(\matrix QR'-RQ'\\-PR'+RP'\\PQ'-QP'\endmatrix\right).
 \label{eq:t}
 \end{equation}

 \end{enumerate}
 \end{LEM}

 \begin{pf} (1) An almost obvious calculation: because $\sI_{C,Z}=(x,y,z)$,
there is a surjective map $\fie\colon\Oh_C^{\oplus3}\to\sI_{C,Z}$, such that
$(h_1,h_2,h_3)\in\ker\fie$ if and only if $h_1x+h_2y+h_3z=0\in\Oh_C$. Write
$H_1,H_2,H_3\in\Oh_\aff$ for lifts of the $h_i$. Then
$H_1x+H_2y+H_3z\in\sI_{\aff^3,C}=(F,G)$. Subtracting off multiples of $F$ and
$G$ means exactly subtracting multiples of the first two rows of $M$ from
$(H_1,H_2,H_3)$, to give identities $H_1'x+H_2'y+H_3'z=0\in\Oh_\aff$. Now
$x,y,z\in\Oh_C$ is a regular sequence, so it follows that $(H'_1,H'_2,H'_3)$ is
in the image of the Koszul matrix given by the bottom 3 rows of $M$. This
proves (1).

(2) A homomorphism $s\colon\sI_{C,Z}\to\Oh_C$ is determined by
$(x,y,z)\mapsto(a,b,c)$ where $a,b,c\in\Oh_C$ satisfy $M(a,b,c)^\mathrm{tr}=0$
(we write $(a,b,c)^\mathrm{tr}$ for the column vector). It is easy to check
that (\ref{eq:t}) gives a map $t$ in this way.

The condition $M(a,b,c)^\mathrm{tr}=0$ consists of 5 equalities in
$\Oh_C=\Oh_{\aff}/(F,G)$. We choose lifts $A,B,C$ to $\Oh_{\aff}$, and write
out the last 3 of these as identities in $\Oh_{\aff}$:
 \begin{equation}
 \begin{array}{cccl}
&-zB&+yC&=\al F -\al'G \\
zA&&-xC&=\be F -\be'G \\
-yA&+xB&&=\ga F -\ga'G
 \end{array}
\quad\text{for some $\al,\dots,\ga'\in\Oh_{\aff}$.}
 \label{eq:al}
 \end{equation}
Taking $x$ times the first plus $y$ times the second plus $z$ times the
third, the left-hand sides cancel, giving the identity
 \begin{equation}
(\al x+\be y+\ga z)F=(\al'x+\be'y+\ga'z)G\in\Oh_{\aff}.
 \nonumber
 \end{equation}
Now since $F,G$ is a regular sequence in $\Oh_{\aff}$, this implies that
 \begin{equation}
 \begin{aligned}
\al x+\be y+\ga z&=DG=D(P'x+Q'y+R'z)\\
\al'x+\be'y+\ga'z&=DF=D(Px+Qy+Rz)
 \end{aligned}
 \label{eq:D}
 \end{equation}
for some $D\in\Oh_{\aff}$.

Now subtracting $D$ times the given generator $t$ changes
 \begin{equation}
\left(\begin{matrix}
A\\B\\C
\end{matrix}\right)
\mapsto
\left(\begin{matrix}
A\\B\\C
\end{matrix}\right)
-\left(\begin{matrix}
QR'-Q'R\\-PR'+P'R\\PQ'-P'Q\\
\end{matrix}\right)D
 \nonumber
 \end{equation}
and has the following effect on the quantities $\al,\dots,\ga'$ introduced in
(\ref{eq:al}):
 \begin{align*}
 (\al,\be,\ga)&\mapsto(\al+DP',\be+DQ',\ga+DR'),\\
 (\al',\be',\ga')&\mapsto(\al'+DP,\be'+DQ,\ga'+DR).
 \end{align*}
To see this, note that the first equation of (\ref{eq:al}) is
 \begin{equation}
-zB+yC=\al F -\al'G=\al(Px+Qy+Rz)-\al'(P'x+Q'y+R'z),
 \nonumber
 \end{equation}
so that the effect of the two substitutions $\al\mapsto\al+DP'$ and
$\al'\mapsto\al'+DP$ on the right exactly cancels out $B\mapsto B+D(PR'-P'R)$
and $C\mapsto C-D(PQ'-P'Q)$ on the left. The upshot is that we can assume
$D=0$ in (\ref{eq:D}).

But then since $(x,y,z)$ is a regular sequence, (\ref{eq:D}) with $D=0$ gives
 \begin{equation}
 \begin{array}{rccc}
 \al=&&ly&-mz\\
 \be=&-lx&&+nz\\
 \ga=&mx&-ny
 \end{array}
\quad\text{and}\quad
 \begin{array}{rccc}
 \al'=&&l'y&-m'z\\
 \be'=&-l'x&&+n'z\\
 \ga'=&m'x&-n'y
 \end{array}
 \nonumber
 \end{equation}
for some $l,\dots,n'\in\Oh_\aff$.
Finally (\ref{eq:al}) can now be rearranged as

 \begin{equation}
 \begin{aligned}
(C-lF+l'G)y&=(B-mG+m'G)z\\
(C-lF+l'G)x&=(A-nF+n'G)z\\
(B-mF+m'G)x&=(A-nF+n'G)y
 \end{aligned}
\quad\text{therefore}\quad
 \begin{aligned}
A-nF+n'G&=Ex\\
B-mF+m'G&=Ey\\
C-lF+l'G&=Ez
 \end{aligned}
 \nonumber
 \end{equation}
for some $E\in\Oh_\aff$. This means that the map $s$ given by $(a,b,c)$ is a
linear combination of $t$ and the identity, as required. \QED \end{pf}

A less pedestrian method of arguing is to say that all three of $\Oh_\aff$,
$\Oh_C$ and $\Oh_Z$ are Gorenstein, so that adjunction gives
 \begin{equation}
0\to\om_C\to\sHom(\sI_Z,\om_C)\to\om_Z=\Ext^1(\Oh_Z,\om_C)\to0. 
 \nonumber
 \end{equation}
The two generators $\id$ and $t$ correspond naturally to the generators of
$\om_C$ and $\om_Z$.

 \pagebreak
 \setlength\parindent{0pt}

Fabrizio Catanese,\par
Dipartimento di Matematica,\par
via Buonarroti 2, I--56127 Pisa (Italy)\par
{\em E-mail address:} catanese@@dm.unipi.it

 \medskip
Marco Franciosi,\par
Scuola Normale Superiore,\par
piazza dei Cavalieri 7, I--56126 Pisa (Italy)\par
{\em E-mail address:} francios@@cibs.sns.it, francios@@gauss.dm.unipi.it

 \medskip
Klaus Hulek,\par
Institut f\"ur Mathematik, Univ.\ Hannover,\par
Postfach 6009 D--30060 Hannover (Germany)\par
{\em E-mail address:} hulek@@math.uni-hannover.de

 \medskip
Miles Reid,\par
Math Inst., Univ.\ of Warwick,\par
Coventry CV4 7AL (England)\par
{\em E-mail address:} Miles@@Maths.Warwick.Ac.UK

 \medskip
Old Uncle Tom Cobbley,\par
Dept.\ of Veterinary Spectrology, Grey Mare Univ.,\par
Widecombe-in-the-Moor, EX31 2KX (England)\par
{\em E-mail address:} Cobblers@@Nag.Agr.UK

 \end{document}